\documentclass[12pt]{article}

\usepackage{amsmath,amsfonts}
\usepackage{amssymb} 
 \usepackage[arrow,matrix,curve]{xy}
\def\xyma{\xymatrix@M.7em}
 \newdir{>>}{{}*!/3.5pt/:(1,-.2)@^{>}*!/3.5pt/:(1,+.2)@_{>}*!/7pt/:(1,-.2)@^{>}*!/7pt/:(1,+.2)@_{>}}
\newdir{ >>}{{}*!/8pt/@{|}*!/3.5pt/:(1,-.2)@^{>}*!/3.5pt/:(1,+.2)@_{>}}
\newdir{ |>}{{}*!/-3.5pt/@{|}*!/-8pt/:(1,-.2)@^{>}*!/-8pt/:(1,+.2)@_{>}}
\newdir{ >}{{}*!/-8pt/@{>}}
\newdir{>}{{}*:(1,-.2)@^{>}*:(1,+.2)@_{>}}
\newdir{<}{{}*:(1,+.2)@^{<}*:(1,-.2)@_{<}}

\newtheorem{satz}{Theorem}[section]

\newtheorem{kor}[satz]{Corollary}
\newtheorem{lem}[satz]{Lemma}

\newtheorem{bem}[satz]{Remark}

\newtheorem{prop}[satz]{Proposition}

\newtheorem{conv}[satz]{Convention} 
\newtheorem{rem}[satz]{Remark}

\newcounter{Roma}
\newenvironment{rom}{\begin{list}{(\roman{Roma})\hfill}{\usecounter{Roma}
\labelsep3mm \leftmargin1.5cm \labelwidth8mm}}{\end{list}}

\makeindex

\newcommand{\nc}{\newcommand}

\nc{\mapco}{\,\colon\, }
\nc{\ab}{^{ ab}}

\nc{\comment}[1]{
}

\nc{\catc}{{\C}}
\nc{\we}{\vee}

\nc{\hra}{\hookrightarrow}

\nc{\epi}{epimorphism}
\nc{\repi}{regular epimorphism}
\nc{\mono}{monomorphism}
\nc{\iso}{isomorphism}
\nc{\coker}[1]{\mbox{${\rm Coker}(#1)$}}
\nc{\Ker}[1]{\mbox{{\rm Ker}$(#1)$}}
\nc{\defgl}{\stackrel{def}{=}}

\nc{\V}{\vspace{3mm}}
\nc{\VV}{\vspace{4mm}}
\nc{\lra}{\longrightarrow}
\nc{\lla}{\longleftarrow}
\nc{\mr}[1]{ \stackrel{#1}{\lra} }
\nc{\ml}[1]{ \stackrel{#1}{\lla} }
\nc{\hmr}[1]{\hspace{2mm}\stackrel{#1}{\lra}\hspace{2mm}}
\nc{\hml}[1]{\hspace{2mm} \stackrel{#1}{\lla}\hspace{2mm}}
\nc{\N}{\noindent}
\nc{\st}{^{\prime}}
\nc{\ot}{\otimes}
\nc{\hcong}{ \hspace{2mm}\cong\hspace{2mm}  }
\nc{\hfbox}{\hfill$\Box$}
\nc{\REF}[1]{(\ref{#1})}

\nc{\Z}{\mbox{$\mathbb{Z}$}} 
\nc{\Q}{\mbox{$\mathbb{Q}$}}

\def\Q{\ifmmode{Q\hskip-5.0pt\vrule height6.0pt depth
0pt\hskip6pt}\else{\hbox{$Q\hskip-5.0pt\vrule height6.0pt depth
0pt\hskip6pt$}}\fi}

\nc{\Ph}{\phantom{}}

\nc{\BE}{\begin{equation}}
\nc{\EE}{\end{equation}}

\nc{\dst}{\displaystyle}
\nc{\sst}{\scriptscriptstyle}
\nc{\ssst}{\scriptscriptstyle}
\nc{\proof}{\N{\bf Proof\,:}\quad}
\nc{\proofof}[1]{\N{\bf Proof of {#1}\,:}\quad}
\nc{\proofofthm}[1]{\N{\bf Proof of Theorem \ref{#1}\,:}\quad}

\nc{\htt}[1]{^{\otimes #1}}

\newcommand{\ssur}[2]{\mbox{$#1 \!\to\!\!\!\!\!\to\! #2$}}

\nc{\Sur}[1]{\mbox{$\:\stackrel{#1}{\lra\!\!\!\!\!\to\,}\:$}}

\def\INJ{\mbox{\mathsurround=0pt
\makebox[0mm][r]{\parbox{0mm}{\rule[-0.65mm]{0mm}{0.2mm}$\scriptscriptstyle>$}}
\makebox[0.7cm][l]{\parbox{0.7cm} {$\lra$}}}}

\def\Inj#1{\mbox{$\:\stackrel{#1}{\INJ}\:$}}

\nc{\Injup}[1]{\mapup{#1}}
\nc{\injup}[1]{\mapup{#1}}
\nc{\injdown}[1]{\mapdown{#1}}

\def\mapup#1{\mbox{$\rule[-1mm]{0cm}{0.7cm}  
\makebox[0mm][r]{\raisebox{0.2mm}{$\scriptstyle\phantom{\cong}$}\hspace{0.6mm}}
\bigg\uparrow\rlap{$\vcenter{\hbox{$\scriptstyle#1$}}$}$}}

\def\isoup#1{\mbox{$\rule[-1mm]{0cm}{0.7cm}  
\makebox[0mm][r]{\raisebox{0.2mm}{$\scriptstyle\cong$}\hspace{0.6mm}}
\bigg\uparrow\rlap{$\vcenter{\hbox{$\scriptstyle#1$}}$}$}}

\def\mapdown#1{\mbox{$\rule[-1mm]{0cm}{0.7cm}  
\makebox[0mm][r]{\raisebox{0.2mm}{$\scriptstyle $}\hspace{0.6mm}}
\bigg\downarrow\rlap{$\vcenter{\hbox{$\scriptstyle#1$}}$}$}}

\def\isor#1{\mbox{$\smash{\mathop{\longrightarrow}\limits^{\cong}_{#1}}$}}

\def\surdown#1{\makebox[0mm]{$\mapdown{\raisebox{0.4mm}{$\scriptstyle#1$}}$}
\makebox[0mm]{\raisebox{-2.15mm}{$\downarrow$}}}

\def\surup#1{\makebox[0mm]{$\mapup{\raisebox{0.4mm}{$\scriptstyle#1$}}$}
\makebox[0mm]{\raisebox{0.7mm}{$\mapup{}$}}}

\newcommand{\brokenr}[2]{
@\cdhgeneric>\raise2.8pt\hbox to3.5pt{\hrulefill}\mkern9mu>
\raise2.8pt\hbox to3.5pt{\hrulefill}\hbox to5pt{}>
\mkern-7.5mu\dashrightarrow>#1 >
#2>}

\newcommand{\brokenup}[1]{
@\cdvgeneric>\hat\cdot
>\raisebox{3pt}{$\vdots$}>
\vbox{\kern3.5pt\hbox{$\cdot$}\kern-3.5pt}
> >\hspace{1mm}#1>}

\newcommand{\functlr}[2]{
\raisebox{0.4pt}{$\hss\begin{CD}
@>\vbox{\hbox to 0pt{$\hss\begin{CD}@<#1<<\end{CD}\hss$}\vskip-2pt}
>#2 >
\end{CD}\hss$}
}

\newcommand{\maprr}[2]{
\raisebox{-0.9pt}{$\hss\begin{CD}
@>\vbox{\hbox to 0pt{$\hss\begin{CD}@>#1>>\end{CD}\hss$}\vskip-3pt}
>#2 >
\end{CD}\hss$}
}

\newcommand{\mapdd}[2]{
@\cdvstandard>\downarrow\hspace{1.5pt}\kern-4pt\hspace{1.5pt}\Big\downarrow>#1>#2>}

\newcommand{\mapud}[2]{
@\cdvstandard>\uparrow\kern-3pt\Big\downarrow>#1>#2>}

\newcommand{\sepi}[3]{\,\mbox{$#1\,: \ssur{#2}{#3}$}\,}

\nc{\auf}{\twoheadrightarrow}
 
\nc{\ruled}{\rule[-4mm]{0mm}{0mm}}

\nc{\qu}{quadratic}
\nc{\GBoGB}{G/BG\st \otimes G/BG\st}
\nc{\lstar}{_{\raisebox{-1mm}{$*$}}}
\nc{\sm}{\:{ \wedge}\:}

\nc{\Rmod}{${\bf R}$-module}

\nc{\map}[3]{\mbox{$#1 \mapco #2 \to #3$}}

\nc{\rond}{{\,\sst \circ\,}}

\nc{\ruleu}{\rule{0mm}{7mm}}

\nc{\T}[1]{\tilde{#1}}

\nc{\Imm}[1]{\mbox{${\rm Im}(#1)$}}

\nc{\tw}{\end{document}}

\nc{\QJG}{\frac{\dst I(G) J}{\dst I^2(G) J}}
\nc{\QJH}{\frac{\dst I(H) J}{\dst I^2(H) J}}
\nc{\otz}{\ot}
\nc{\UL}[2]{{\rm U}_{#1}{\rm L}(#2)}

\nc{\IRN}[1]{I_{R,\cal G}^{#1}(G)}

\nc{\ULG}[1]{{\rm U}_{#1}{\rm L}^{\cal G}(G)}
\nc{\UGH}[1]{{\rm U}_{#1}^{\cal GH}(G,H)}

\nc{\calG}{{\cal G}}
\nc{\AB}{^{ AB}}

\nc{\IZ}{I}

 \nc{\UG}{\mbox{${\rm UL}^{\cal G}(G)$}}

 \nc{\UN}{\mbox{${\rm UL}^{\cal N}(N)$}}

 \nc{\UT}{\mbox{${\rm UL}^{\gamma}(T)$}}

 \nc{\UGn}[1]{\mbox{${\rm U}_{#1}{\rm L}^{\cal G}(G)$}}

 \nc{\UNn}[1]{\mbox{${\rm U}_{#1}{\rm L}^{\cal N}(N)$}}

 \nc{\UTn}[1]{\mbox{${\rm U}_{#1}{\rm L}^{\gamma}(T)$}}

\nc{\LG}{\mbox{$ {\rm L}^{ \cal G}(G)$}}
\nc{\LH}{\mbox{$ {\rm L}^{ \cal H}(H)$}}
\nc{\LN}{\mbox{$ {\rm L}^{ \cal N}(N)$}}

\nc{\LnG}[1]{\mbox{$ {\rm L}_{#1}^{\sst\cal G}(G)$}}

\nc{\LNn}[1]{\mbox{$ {\rm L}_{#1}^{\sst\cal N}(N)$}}
\nc{\LnN}[1]{\mbox{$ {\rm L}_{#1}^{\sst\cal N}(N)$}}

\nc{\KnD}[1]{{\cal K}_{#1}^{\Delta}}

\nc{\KnL}[1]{{\cal K}_{#1}^{\Lambda}}

\nc{\KnINL}[1]{{\cal K}_{#1}^{I(N)\Lambda}}

 \nc{\QnG}[1]{Q_{#1}(G)}

 \nc{\QnGT}[1]{Q_{#1}(G,T)}

 \nc{\QnGN}[1]{Q_{#1}(G,N)}

\nc{\Tor}[2]{{\rm Tor}_1^{\mathbb{Z}}(#1\,,#2)}

\def\isor#1{\mbox{$\smash{\mathop{\longrightarrow}\limits^{\cong}_{#1}}$}}

\nc{\caln}{{\cal N}}
\nc{\UNN}[1]{{\rm U}^{\cal N}_{#1}(N,N)}
\nc{\ULN}[1]{{\rm U}_{#1}{\rm L}^{\cal N}(N)}

\nc{\calH}{{\cal H}}
\nc{\calK}{{\cal K}}

\nc{\ULH}[1]{{\rm U}_{#1}{\rm L}^{\cal H}(H)}

\nc{\hoplus}{\hspace{3mm}\oplus\hspace{3mm}}

\begin{document} 

\title{On Fox and augmentation quotients of semidirect products}
\author{Manfred Hartl}
\maketitle

\begin{center}

\N \rm Universit\'e de Lille  Nord de France, 59000 Lille, France\\
UVHC, LAMAV  and FR CNRS 2956,\\ 
59313 Valenciennes, France.\\ 
Email: Manfred.Hartl@univ-valenciennes.fr\\ 

\end{center}
\vspace{8mm}

\begin{abstract}

\N Let $G$ be a group which is the semidirect product of a normal subgroup $N$ and some subgroup $T$.  
Let $I^n(G)$, $n\ge 1$, denote the powers of the augmentation ideal $I(G)$ of the group ring
$\Z(G)$. Using homological methods the groups
$Q_n(G,H) = I^{n-1}(G)I(H)/I^{n}(G)I(H)$, $H=G,N,T$, are functorially expressed in terms of enveloping algebras of certain Lie rings associated with $N$ and $T$, in the following cases: for $n\le 4$ and arbitrary $G,N,T$ (except from one direct summand of $Q_4(G,N)$), and for all $n\ge 2$ if certain filtration quotients of $N$ and $T$ are torsion-free.

\end{abstract}\vspace{5mm}

 \N{\large \bf Introduction.}\quad The group ring $\Z(G)$ of a group $G$ is naturally filtered by the powers $I^n(G)$, $n\ge 1$, of its augmentation ideal $I(G)$. It is a long-studied problem to determine the so-called augmentation quotients $Q_n(G)= I^n(G)/I^{n+1}(G)$ in terms of the structure of $G$, also because of their close link with the dimension subgroups $D_n(G) = G \cap (1+I^n(G))$ which can be inductively described as $D_{n+1}(G) = {\rm Ker}(D_n(G) \to Q_n(G))$. The groups $Q_n(G)$ were determined for $n=2$ by Passi \cite{PaPF}, Sandling \cite{SaR} and Losey \cite{Lo} for abelian, finite and finitely generated groups $G$ and for $n=3,4$ and finite $G$ by Tahara \cite{Ta3}, \cite{Ta4}; a functorial description for all groups was given for $n=2$ by Bachmann and Gruenenfelder \cite{Ba-GrII} and  for $n=3$ in \cite{Q3}, based on Quillen's approximation of the graded ring Gr$(\Z(G)) = \Z \oplus \bigoplus_{n\ge 1} Q_n(G)$  by  the enveloping ring of the Lie ring of $G$, see \cite{Qu}, \cite{Pa} or section 1 below. 
 With the functorial viewpoint new methods emerge in the field where most work is still based on the classical combinatorial approach: Quillen's computation of Gr$(\Q(G))$ relies on Hopf algebra techniques, while the general results on $Q_2(G)$ and $Q_3(G)$ are obtained by homological methods as inaugurated by Passi, see \cite{Pa}. Recently, Passi and Mikhailov pass from homological to homotopical methods \cite{MP}, \cite{MPbook}.
  In \cite{Fox2} and in this paper we combine enveloping algebras and new homological observations to study the more general quotients $Q_n(G,H) = I^{n-1}(G)I(H)/I^{n}(G)I(H)$ for some subgroup $H$ of $G$; we call these {\em Fox quotients}\/ because of their close relation with the classical {\em Fox subgroups}\/ $G \cap  (1+I^{n-1}(G)I(H))$.
  Fox quotients (and some related groups, see \cite{Kh14}, \cite{KV27}, \cite{KV29}) were also extensively  studied in the literature, but, except from  \cite{Fox2}, only under suitable splitting assumptions, in particular when $H$ is a semidirect factor of $G$. In fact, Sandling's \cite{Sa} and later Tahara's work \cite{Ta} on augmentation quotients of semidirect products $G=N\rtimes T$ had split the study of Fox quotients into  two classes of independent problems: the study of certain filtration quotients of  $\Z(N)$ and $\Z(T)$ on the one hand and of product filtrations ${\cal F}_n=\sum \Delta_{n-i} I^i(T)$ on the other hand where $(\Delta_i)_{i\ge 1}$ is one of two natural filtrations of $\Z(N)$, see section 1. In a series of papers   Khambadkone and later Karan and Vermani expressed the quotients of these product filtrations in terms of tensor products of the groups $\Delta_{n-i}/\Delta_{n-i+1}$ and $I^i(T)/I^{i+1}(T)$, for low values of $n$ and under additional assumptions, assuming either $G$ finite and $N$ finitely generated or nilpotent \cite{Kh13}, \cite{Kh14}, \cite{Kh15}, or assuming torsion-freeness of sufficiently many filtration quotients of $N$ and $T$ \cite{KV30}, \cite{KV32}, \cite{KV31},  \cite{KV47}. For a more detailed survey on  Fox and augmentation quotients  see Passi \cite{Pa} and Vermani \cite{Ve}.
 
 In this paper we treat the general case, showing that the quotients of the product filtrations above are in fact iterated amalgamated sums of tensor products 
 of the groups $\Delta_{n-i}/\Delta_{n-i+1}$ and $I^i(T)/I^{i+1}(T)$, amalgamated along certain subgroups of torsion products of these groups. We thus completely determine the groups $Q_n(G,H)$ for arbitrary $H=G,N,T$ and $n\le 3$, and for $n=4$ with the exception of one of the two direct factors of $Q_4(G,N)$, see section 2. Our description is functorial and in terms of enveloping rings of certain Lie algebras associated with $N$ and $T$, see section 1. If suitable filtration quotients of $N$ and $T$ are torsion-free then our amalgamated sums degenerate to direct sum decompositions; we then express $Q_n(G,H)$ for $H=G,N,T$ 
 in terms of tensor products of enveloping rings as above, thus improving and generalizing similar results of Karan and Vermani for $n\le 4$ to all $n\ge 1$, see section 3.
 
  The first three sections are devoted to a presentation of the necessary constructions and results while the proofs are deferred to sections 4 and 5. \V

\section{Preliminary constructions and results}

In this section we recall and partially generalize constructions and results from the literature which are needed in the sequel.

Let $G$ be a group. An {\em N-series}\/ ${\cal G}$  of $G$
is a descending chain of subgroups
  \[   G= G_{(1)} \supset G_{(2)}\supset G_{(3)}\supset \ldots  \]
such that $[ G_{(i)}, G_{(j)}] \subset G_{(i+j)}$ for $i,j\ge 1$, with $[a,b] = aba^{-1}b^{-1}$.
A given  N-series ${\cal G}$ induces a descending chain of two-sided ideals of the group ring $\Z(G)$
  \[  \Z(G) = I^0_{ \cal G}(G) \supset I^1_{ \cal G}(G) \supset I^2_{ \cal G}(G) \supset \ldots  \]
by defining $I^n_{ \cal G}(G)$  (for $n\ge 1$) to be the  subgroup of $\Z(G)$ generated by the elements
  \[  (a_1 -1) \cdots (a_r -1)\,, \quad \mbox{$r\ge 1$, $a_i \in G_{(k_i)}$, such that $k_1+ \ldots +
k_r  \ge n$.}  \] 
Two examples of N-series are used  throughout in this paper:
\begin{itemize}
\item  the lower central series $\gamma=(G_i)_{i\ge 1}$, $G_1=G$ and $G_{i+1}=[G_i,G]$,  where the inclusion
$I(G_n)\subset I^n(G)$ implies that
$I^n_{\gamma}(G)$ equals $I^n(G)$, the $n$-th power of the ideal $I(G)$;

\item for a normal subgroup $N$  of $G$  an N-series ${\cal N} = (N_{(i)})_{i\ge 1}$ of $N$
is defined by $N_{(i)} = [N_{(i-1)},G]$; note that if $G=NT$ for some subgroup $T$ then $N_{(i)} = [N_{(i-1)},N][N_{(i-1)},T]$.

\end{itemize}

The second
example was introduced by Tahara; we here gather some basic results due to him \cite{Ta} and to Khambadkone \cite{Kh13}.\vspace{2mm}


\begin{satz}\label{prelimsplit} Suppose that $G$ is the semidirect product of a normal subgroup $N$ and some subgroup $T$. Write
$\Lambda_n = I^n_{\cal N}(N)$ and let 
  \[ {\cal K}_n = \sum_{i=1}^{n-1} \Lambda_{n-i} \, I^i(T)\:, \hspace{8mm} {\cal K}_n^* = \sum_{i=1}^{n} \Lambda_{n-i} \, I^i(T)\:,
   \hspace{8mm} \Gamma_n^* = \sum_{i=0}^{n-1}
I^{n-i}(T)  \, \Lambda_{i} \:.\]
Then   
  \BE\label{NT} I(G) = I(N) \hoplus I(T) \hoplus I(N)I(T) \EE
  \BE\label{TN} \phantom{I(G)} = I(T) \hoplus I(N) \hoplus I(T)I(N) \EE
  \BE\label{AugNT} I^n(G) = \Lambda_n \hoplus I^n(T) \hoplus {\cal K}_n \EE
  \BE\label{GT} I^n(G) I(T) =  I^{n+1}(T) \hoplus {\cal K}_{n+1} \EE
  \BE\label{GN} I^n(G) I(N) = \Lambda_n I(N) \hoplus   \Gamma_n^*I(N) \EE
  \BE\label{NG} I(N) I^n(G)  = I(N) \Lambda_n  \hoplus   I(N)  {\cal K}_n^*\:. \EE
\end{satz}

\comment{ 
\[
\begin{array}{ccccccc}
   \label{NT} I(G) &=& I(N) & \hoplus & I(T)  & \hoplus &  I(N)I(T) \\
   \label{TN}      &=& I(T) & \hoplus & I(N)  & \hoplus &  I(T)I(N) \\
   \label{AugNT} I^n(G) &=& \Lambda_n  & \hoplus &  I^n(T)  & \hoplus &  {\cal K}_n \\
   \label{GT} I^n(G) I(T) &=&   &   & I^{n+1}(T)  & \hoplus &  {\cal K}_{n+1} \\
   \label{GN} I^n(G) I(N) &=&  & & \Lambda_n I(N)  & \hoplus &   \Gamma_n^*I(N)\:. 
\end{array}
\]
}

We need to make the obvious symmetry between the relations \REF{GN} and \REF{NG} more precise.

\begin{rem}\label{conjonFox}\rm Let $G$ be a group with  distinguished N-series $\calG$ and subgroup  $H$. Then
the anti-ring-automorphism $(-)^{\star}\colon \Z(G) \to \Z(G)$ sending $g\in G$ to $g^{-1}$, called conjugation, carries $I(H)I^n_{\calG}(G)$ onto 
$I^n_{\calG}(G)I(H)$ since the subgroups $H$ and $G_{(i)}$ are stable under inversion. Hence it induces an isomorphism of abelian groups
\[ I(H)I^{n-1}_{\calG}(G)/I(H)I^n_{\calG}(G) \hcong
I^{n-1}_{\calG}(G)I(H)/ I^n_{\calG}(G)I(H)\,. \]
For $G=N\rtimes T$ as above, we also have 
\begin{eqnarray}
(I(N)  {\cal K}_n^* )^{\star} &=& \Big(I(N) \sum_{i=1}^{n} \Lambda_{n-i} \, I^i(T)\Big)^{\star} \nonumber\\
&=&\sum_{i=1}^{n}I^i(T)\,\Lambda_{n-i}I(N)\nonumber\\
&=  &
  \sum_{j=0}^{n-1} I^{n-j}(T)\,\Lambda_{j}I(N)\nonumber\\
&=  &\Gamma_n^*I(N)  \,. \label{conjKappatoGamma}
\end{eqnarray}

\end{rem}

 Now recall that our aim is to determine the filtration quotients
\[Q_n(G ) = I^{n }(G) /I^{n+1}(G) \]\[
 Q_n(G,H) = I^{n-1}(G)I(H)/I^{n}(G)I(H)  \]
for $H=N,T$; note that $Q_n(G) = Q_n(G,G)$.
The relations above immediately imply the following identities.
  \BE\label{AugTallg} Q_n(G) \hspace{3mm}=\hspace{3mm}  \Lambda_n/\Lambda_{n+1} \hoplus Q_n(T) \hoplus {\cal K}_n/{\cal K}_{n+1} \EE
  \BE\label{FoxTallg} Q_n(G,T) \hspace{3mm}=\hspace{3mm}    Q_n(T) \hoplus {\cal K}_n/{\cal K}_{n+1} \EE
  \BE\label{FoxNallg} Q_n(G,N) \hspace{3mm}=\hspace{3mm}  \Lambda_{n-1}I(N)/\Lambda_{n}I(N) \hoplus  
\Gamma_{n-1}^*I(N)/\Gamma_n^*I(N) \ruled \EE 
It turns out that the terms on the right hand side of the above identities fall into two categories:
first of all,   ${\cal K}_n/{\cal K}_{n+1}$ and 
$\Gamma_{n-1}^*I(N)/\Gamma_n^*I(N)$ each of which arises from the product  of a filtration  of $I(N)$ with one of $I(T)$; the strategy here is to express the quotients of these product filtrations in terms of - tensor and torsion - products of the factors. Once this "separation of the factors" is achieved (which is the main concern of this paper, see sections 2 and 3) one is left with dealing with the generalized Fox and augmentation quotients 
of $N$ and $T$, i.e. the groups  $Q_n^{\calK}(K) = I^n_{\calK}(K)/I^{n+1}_{\calK}(K)$ for $K=N,T$  and $\calK={\cal N},\gamma$, resp., and
$Q_n^{\cal N}(N,K) = I^{n-1}_{\cal N}(N)I(K)/I^{n}_{\cal N}(N)I(K)$ for    some subgroup $K$ of $N$ (here only the case $K=N$ is needed), see the results in sections 2 and 3.  The study of these groups requires the following constructions.

The basic idea, due to Quillen \cite{Qu}, is to approximate the groups $Q_n^{\calG}(G)$ by means of enveloping algebras. The
construction for arbitrary N-series $\calG$ can be found in Passi's book \cite{Pa}, but we recall it here for convenience of the
reader and to fix notation.

The abelian group ${\rm L}^{\cal G}(G) = \sum_{i\ge 1} G_{(i)}/G_{(i+1)}$ is a graded Lie ring whose bracket is induced by the
commutator pairing of
$G$. So its enveloping algebra $\ULG{}$ over the integers is defined. On the other hand, the  filtration quotients
$Q_n^{\calG}(G) = I^n_{\cal G}(G)/I^{n+1}_{\cal G}(G)$ form the graded ring  Gr$^{\cal G}(\Z(G)) = \oplus_{i=0} Q_n^{\calG}(G)$;
note that one has  Gr$^{\gamma}(\Z(G))=  \bigoplus_{n\ge 0} I^n(G)/I^{n+1}(G)$. \comment{If $\calG =\gamma$ we frequently
suppress the superscript $\gamma$ from notation.} Now the map 
$\LG{} \to$ Gr$^{\cal G}(\Z(G))$, $aG_{(i+1)} \mapsto a-1+I^{i+1}_{\calG}(G)$ for $a\in G_{(i)}$, is a homomorphism of graded Lie
rings and hence extends to a map of graded rings 
  \[ \theta^{\cal G}\,\colon\, \UG \Sur{} {\rm Gr}^{\cal G}(\Z(G)) \:.\]
 This map is clearly
surjective but rarely globally injective; for instance, $\theta^{\gamma}$ is  
injective  if $G$ is cyclic, but is non
injective for all non cyclic finite abelian groups \cite{Ba-GrI}. An important favourable case is given by the following result which
relies on work of Hartley \cite{Hartley}.\vspace{2mm}

\begin{satz}\label{theta} \quad Let $G$ be a group and let ${\cal G}\,\colon\, G = G_{(1)}
\supset G_{(2)} \supset \cdots$ be an N-series of $G$ with torsion-free quotients $ G_{(i)}/
G_{(i+1)}$ for all $i\ge 1$. Then $\theta^{\cal G}$ is an isomorphism.
\end{satz}

\proof In a first step, one adapts an argument of Quillen \cite{Qu} for the case ${\cal
G} = \gamma$ to show that the epimorphism
   \[ \theta^{\cal G}_{Q\!\!\!\!I} = \theta^{\cal G}\ot \Q\,\colon\, {\rm U}( \LG \ot \Q) \Sur{} {\rm Gr}^{\cal G}(\Q(G)) \]
is an isomorphism: indeed, the image of the canonical map of graded Lie rings \map{j_{\cal G}}{\LG}{{\rm Gr}^{\cal
G}(\Q(G))} consists of primitive elements (with respect to the Hopf algebra structure induced by
the canonical one of $\Q(G)$), and generates ${\rm Gr}^{\cal G}(\Q(G))$  as an algebra (by
definition of the filtration $(I_{\Q,\cal G}^i(G))_{i\ge 0}$). So by the Milnor-Moore theorem
${\rm Gr}^{\cal G}(\Q(G))$ can be identified with the enveloping algebra of the Lie algebra of its primitive elements.
Hence it suffices to check that the map $ j_{\cal G} \ot \Q$ is injective since then by a standard argument the Poincar\'e-Birkhoff-Witt theorem implies
that $\theta^{\cal G}\ot \Q$ is injective, too: 
to recall this, let $f\colon {\frak g}\to {\frak g}'$ be an injective map of Lie algebras over some field $\mathbb{K}$; we wish to show that the map of $\mathbb{K}$-algebras ${\rm U}(f)\colon {\rm U}({\frak g}) \to {\rm U}({\frak g}')$ induced by $f$ is also injective. Denote by $ {\rm SP}({\frak g})$ the symmetric algebra over the $\mathbb{K}$-vector space ${\frak g}$, and by ${\rm GR}({\rm U}({\frak g}))$ the graded algebra associated with the canonical increasing filtration of ${\rm U}({\frak g})$, see \cite{PBW}. Now under the natural Poincar\'e-Birkhoff-Witt  isomorphism of graded algebras $ {\rm SP}({\frak g}) \stackrel{\cong}{\to}
{\rm GR}({\rm U}({\frak g}))$, the map  ${\rm GR}({\rm U}(f)) \colon {\rm GR}({\rm U}({\frak g})) \to{\rm GR}({\rm U}({\frak g}'))$ corresponds to the  map ${\rm SP}(f)$, which is injective since if $r\colon \frak{g}\to \frak{g}'$ is a $\mathbb{K}$-linear retraction of $f$ then ${\rm SP}(r)$ is a retraction of ${\rm SP}(f)$. Thus ${\rm GR}({\rm U}(f))$
 is injective, too, and hence so is  ${\rm U}(f)$, since any element of ${\rm U}({\frak g})$ lies in a finite step of the increasing filtration.


Now our map
$j_{\cal G}$ is indeed injective by the hypothesis on ${\cal G}$, as is proved in \cite{Hartley}, see also
\cite{Pa}. Thus also $\theta_{\Q}^{\cal G}$ is injective, whence an isomorphism. 

In order to descend to integral coefficients, consider the following sequence of
homomorphisms of graded rings.
 \[  \UG \Sur{\theta^{\cal G}}
{\rm Gr^{\cal G}}({\Z}(G)) \stackrel{\bar{\iota}}{\lra} {\rm Gr}^{\cal G}({\Q}(G))
\,\isor{(\theta_{Q\!\!\!\!I}^{\cal G})^{-1} \ruled }\, {\rm U}(  \LG \otimes \Q)\,
\cong\,  \UG  \otimes \Q \]
Here $\bar{\iota}$ is induced by the canonical injection ${\iota}\colon {\Z}(G) \hra {{\Q}(G)}$.
The composite map $\zeta\colon \UG \to   \UG \ot\Q$ sends $x\in\UG$ to $x\ot 1$; this is easily checked on elements of ${\rm L}^{\cal G}(G)$ which generate \UG\ as a ring. Hence $\zeta$
is injective since $\UG$ is
torsion-free as \LG\ is, see\ \cite[Lemma 1.11]{PolProp}. So the first factor $
\theta^{\cal G}$ of $\zeta$ is also injective, as was to be shown. \hfbox\V

\begin{kor}\label{UGrbism}  Let ${\cal G}$ be an N-series of $G$ such that for some $m\ge 1$ the groups $G_{(i)}/G_{(i+1)}$ are
torsion-free for $1\le i\le m$. Then in this range the map   $\theta_{i}^{\cal G}\,\colon\, \UGn{i}$ $
\Sur{} {\rm Gr}^{\cal G}_i(\Z(G))$ is an isomorphism  and the groups  $I(G)/I^{i+1}_{\cal G}(G)$ are torsion-free.
\end{kor}

\proof  Just note that passing to the quotient $G\auf G/G_{(m+1)}$ does affect neither  \UGn{i}  nor \,Gr$_i^{\cal G}(\Z(G))$ for
$i\le m$; but $G/G_{(m+1)}$ satisfies the hypothesis of Theorem \ref{theta} for the N-series $\pi{\cal G}\,\colon\, G/G_{(m+1)}
\supset  G_{(2)}/G_{(m+1)} \supset\ldots$ of $G/G_{(m+1)}$. Moreover, \,L$^{\pi{\cal G}}(G/G_{(m+1)})$ being torsion-free so is its
  enveloping algebra (see \cite[lemma 1.11]{PolProp}) and hence  
\,Gr$_i^{\cal G}(\Z(G/G_{(m+1)})) \:\cong\:$\,Gr$_i^{\cal G}(\Z(G)) = I^i_{\cal G}(G)/I^{i+1}_{\cal G}(G)$ for $i\le m$.
Consequently also $I(G)/I^{i+1}_{\cal G}(G)$ is torsion-free being an iterated extension of the groups
$I^j_{\cal G}(G)/I^{j+1}_{\cal G}(G)$ for $1\le j\le i$.\hfbox\V

If $\calG$ is an arbitrary N-series, the kernel of $\theta^{\calG}$ is   a torsion group since $\theta^{\cal G}_{Q\!\!\!\!I}$ is an isomorphism,
but its structure remains
widely unknown,  even for $\calG=\gamma$. At least in low degrees the problem is settled; the following result was obtained in \cite{Ba-GrI} for $\calG=\gamma$ and  in \cite{Q3} for arbitrary N-series $\calG$.

\begin{satz}\label{U1-U2} The map $\theta^{\calG}_n$ is an isomorphism for $n=1,2$ and all groups $G$ and N-series $\calG$ of $G$.
\end{satz}

This result is actually generalized in Theorem \ref{thetaGH12} below.

We need to make Theorem \ref{U1-U2} more explicit:  it says that $Q_1^{\calG}(G) \:\cong\: G/G_{(2)}$, and that $Q_2^{\calG}(G)$ can be described as follows.
Recall that the exterior square $A\sm A$ of an abelian group $A$ is defined to be the 
quotient
of $A \ot A$ modulo the subgroup generated by the diagonal elements.
We write  $G\ab=G/G_{2}$ and $G\AB=G/G_{(2)}$, the same for $H$; note that $G\AB=\LnG{1} \cong \ULG{1}$. Then we have natural homomorphisms
\[   G_{(2)}/G_{(3)} \hml{c_2^{\calG}}   G\AB \sm G\AB \hmr{l_2^{\calG}} G\AB \ot G\AB \hmr{\mu_2^{\calG}}
Q_2^{\calG}(G) \hml{p_2^{\calG}} G_{(2)}/G_{(3)} \]
 defined for $a,b \in G$ by $\,c_2^{\calG}(aG_{(2)} \sm b G_{(2)}) = [a,b] G_{(3)}$,   $l_2^{\calG}(aG_{(2)} \sm bG_{(2)}) = aG_{(2)} \ot bG_{(2)} - bG_{(2)} \ot
aG_{(2)}$,
  $\mu_2^{\calG} (aG_{(2)} \ot bG_{(2)}) = (a-1)(b-1) + I^3_{\calG}(G)$, and $p_2^{\calG}(aG_{(3)}) = a-1 + I^3_{\calG}(G)$. 
Then by Theorem \ref{U1-U2} the following sequence is exact: \ruled
  \BE\label{strofU2} G\AB \sm G\AB \hmr{(c_2^{\calG},-l_2^{\calG})} G_{(2)}/G_{(3)} \hoplus G\AB \ot G\AB
\hmr{(p_2^{\calG},\,\mu_2^{\calG})^t} Q_2^{\calG}(G) \hmr{} 0
\EE

In order to describe the third Fox and augmentation quotients of $G$ we need the following constructions. Let $\calG$ be an N-series of $G$,
$H$ a subgroup of $G$, and ${\cal H}=(H_{(i)})_{i\ge 1}$ be an N-series of $H$ such that $H_{(i)} \subset G_{(i)}$ for $i\ge 1$. 
These data give rise to a filtration
\[ {\cal F}^1 = \Z(G)I(H) \supset {\cal F}^2 =\Z(G)I(H_{(2)}) + I(G)I(H) \supset \ldots \]
 of $\Z(G)I(H)$ by sub-$\Z(G)$-$\Z(H)$-bimodules ${\cal F}^n$ defined by 
  \[ {\cal F}^n = \sum_{\begin{matrix} \mbox{$\sst i\ge 0$, $\sst j\ge 1$} \cr \mbox{$\sst i+j=n$} \end{matrix}}
I^{i}_{\cal G}(G) I^{j}_{\cal H}(H)  =  
\sum_{\begin{matrix} \sst \mbox{$\sst i\ge 0$, $\sst j\ge 1$}\cr  \mbox{$\sst i+j=n$} \end{matrix}} I^{i}_{\cal G}(G) I(H_{(j)}) \:.\] 
Note that if ${\cal H} =\gamma$ then ${\cal F}^n = I^{n-1}_{\cal G}(G) I(H)$ since $I^j_{\gamma}(H) = I^j(H)$, and if $H=G$ and $\cal H=G$ then  ${\cal F}^n =
I^{n}_{\cal G}(G)$. 
  The associated graded group \,Gr$^{\cal GH}(\Z(G)I(H)) = \bigoplus_{n\ge 1} {\cal F}^n/{\cal F}^{n+1}$ is  a  graded Gr$^{\cal
G}(\Z(G))$--Gr$^{\calH}(\Z(H))$--bimodule in the canonical way, and hence a $\ULG{}-{\rm UL}^{\calH}(H)$-bimodule via the maps
$\theta^{\calG}$ and $\theta^{\cal H}$.

We now generalize the approximation of the ring\,Gr$^{\calG}(\Z(G))$ by $\ULG{}$ to the bimodule \,Gr$^{\cal GH}(\Z(G)I(H))$, as
follows. The injection $\iota\,\colon\,H \hra G$ induces a canonical map of graded Lie rings ${\rm L}(\iota)\colon \LH \to \LG $ which extends to 
a map of graded rings ${\rm UL}(\iota)\colon \ULH{} \to \ULG{}$. It makes $\ULG{}$ into a $\ULH{}$-bimodule, whence the graded 
 $\ULG{}$-$\ULH{}$-bimodule
\BE\label{UGHdef} \UGH{} =  \ULG{}
\ot_{\raisebox{-1mm}{$\sst {\rm UL}^{\cal H}(H)$}}
\bar{\rm U}\LH \EE
is defined where $\bar{\rm U}\LH$  denotes the augmentation ideal of $\ULH{}$. Now let the surjective map of
$\ULG{}$-$\ULH{}$-bimodules
  \[  \theta^{\cal GH}  \:\colon \hspace{2mm} \UGH{}  
\hspace{2mm}\Sur{}\hspace{2mm} {\rm Gr}^{{\cal GH}}(\Z(G)I(H)) \]
be defined as follows: as $I^n_{\calH}(H) \subset {\cal F}_n$ for $n\ge 1$, we obtain a map of graded left ${\rm Gr}^{\cal H}(\Z(H))$-modules  
${\rm Gr}^{\cal H}(I(H)) \lra  {\rm Gr}^{\cal GH}(\Z(G)I(H))$ which by precomposition with 
$\theta^{\calH}$ gives rise to a map of graded left 
$\ULH{}$-modules $\bar{\rm U}\LH \lra {\rm Gr}^{\cal GH}(\Z(G)I(H))$. Now $\theta^{\cal GH}$ is obtained by extension of scalars along the map ${\rm UL}(\iota)$. More explicitly, for $i\ge 0$, $j\ge 1$ such that $i+j=n$, $x\in \ULG{i}$, $y\in \ULH{j}$, $x\st\in I^{i}_{\cal G}(G)$ and
$y\st \in I^{j}_{\cal H}(H)$ such that $\theta^{\calG}_i(x) = x\st + I^{i+1}_{\cal G}(G)$ and
$\theta^{\calH}_j(y) = y\st + I^{i+1}_{\cal H}(H)$, one has $\theta^{\cal GH}_n(x\ot y) = x\st y\st + {\cal F}^{n+1}$. Note that for
$H=G$ and $ \cal H=G$, $\theta^{\cal GG} = \theta^{\calG} \mu^{\calG}$ where $\mu^{\calG}\,\colon\, {\rm U}^{\cal GG}(G,G)
\mr{\cong} \bar{\rm U}\LG$ is the canonical isomorphism.

We now study the map $\theta^{\cal GH}$ in  degree $n\le 3$. To start with, it is clear that
\BE\label{U1GH} {\rm U}_1^{\cal GH}(G,H) \hspace{6pt}\cong\hspace{6pt}
{\rm U}_1{\rm L}^{\cal  H}(H) \hspace{6pt}\cong\hspace{6pt}{\rm L}_1^{\cal  H}(H) \hspace{6pt}\cong\hspace{6pt} H\AB\,.
\EE

The following convention turns out to be convenient throughout the rest of this paper.

\begin{conv} For a group $K$ with N-series ${\cal K} = (K_{(i)})_{i\ge 1}$ and $a\in K_{(i)}$ we  consider the coset $aK_{(i+1)} \in
K_{(i)}/K_{(i+1)} = {\rm L}_i^{\cal K}(K)$ also as an element of ${\rm U}_i{\rm L}^{\cal K}(K)$, thus suppressing the canonical map 
$\nu\colon {\rm L}^{\cal K}(K) \to {\rm U}{\rm L}^{\cal K}(K)$ from the notation. 
Moreover, all sums, powers and products  of cosets enclosed between brackets, of the type $(aK_{(i+1)})$ with $a\in K_{(i)}$, are understood to be taken in the ring $ {\rm U}{\rm L}^{\cal K}(K)$.

\end{conv}\vspace{2mm}

Now for $i\ge 0$ and $j\ge 1$ let 
   \[  \nu_{ij}\,\colon\,\ULG{i} \ot \ULH{j} \to
{\rm U}_{i+j}^{{\cal GH}}(G,H)  \]
 be the canonical map, and 
 \[ u \colon \bar{\rm U}{\rm L}^{{\calH}}(H) \lra {\rm U}^{{\cal GH}}(G,H)\]
 be the map of graded $\ULH{}$-bimodules given by $u(x) =1\ot x$. Finally, let
\[ q_n^{\calG}\colon G/G_n \Sur{} 
 G/G_{(n)}\]
 be the canonical quotient map induced by the inclusion $G_n \subset G_{(n)}$; and maps denoted by   $\pi_k$ for some or no index $k$ always denote  canonical quotient maps.

 \begin{satz}\label{thetaGH12} The maps $\theta^{\cal G H}_n$ are isomorphisms for $n=1,2$ and all $G,H,\calG$ and compatible $\calH$ as above.
 
 \end{satz}
 
 Note that taking $(H,\calH) =(G,\calG)$ we thus recover Theorem \ref{U1-U2}.\V
 
 \proof  For $n=1$, we use the classical isomorphism
 \BE\label{Fox1} D_H\colon  H\ab\hspace{7pt} \mr{\cong} \hspace{7pt} \frac{\Z(G)I(H)}{I(G)I(H)}\,,\quad D_H(hH_2) = h-1+I(G)I(H)\,,\EE
 due to Whitcomb, see \cite{Whitcomb}, or \cite[Proposition 3.1]{Fox2} for an easy homological proof of a more general fact. We obtain isomorphisms
 \[  {\rm U}_1^{\cal GH}(G,H) \hspace{4pt}\cong\hspace{4pt} H\AB \hspace{4pt}\cong\hspace{4pt}
 \frac{(H/H_2)}{\Imm{H_{(2)}}} \hspace{4pt}\stackrel{\overline{D_H}}{\cong} \hspace{4pt} 
 \frac{\Z(G)I(H)}{I(G)I(H)}\Big/\Imm{\Z(G)I(H_{(2)})} 
 \hspace{4pt}\cong\hspace{4pt} {\cal F}_1/{\cal F}_2\]
 whose composite is $\theta_1^{\cal GH}$.

 Now let $n=2$. As in \cite{Fox2}, our arguments are most conveniently formulated in the language of pushouts of abelian groups, thereby using their elementary properties, in particular the gluing of pushouts and the link between the  kernels of parallel maps in a pushout square, see \cite{Rotman}. 
 
 Consider the following diagram whose top row is exact taking the right-hand map to be given by the  projection to the cokernel of the left-hand map followed by the isomorphism $D_H^{-1}$ in
\REF{Fox1}, and
where $j_1$ is induced by the corresponding injection  and $D'$ is given by restriction of $j_1D_{H_{(3)}}$.
 \BE\label{intersecdia}
 \xymatrix{
\frac{\dst I(G)I(H)}{\dst I^2_{\calG}(G)I(H) + I(G)I^2_{\calH}(H)}\hspace{7pt} \ar@{^(->}[r] & \hspace{7pt}\frac{\dst  \Z(G)I(H)}{\dst I^2_{\calG}(G)I(H) + I(G)I^2_{\calH}(H)}\hspace{7pt} \ar[r] & \hspace{7pt}\frac{\dst H}{\dst H_2}\ar[r]\hspace{7pt} & 0\\
& \frac{\dst  \Z(G)I(H_{(3)})}{\dst I(G)I(H_{(3)})}
\ar[u]^{j_1}& \\
\rule{0mm}{20pt}\rule[-14pt]{0mm}{20pt}\frac{\dst H_{(3)}\cap H_2}{\dst [H_{(3)}, H_{(3)}]} \ar@{^(->}[r] \ar[uu]^{D'}  & \rule{0mm}{20pt}\rule[-14pt]{0mm}{20pt}\rule{0mm}{7mm}\frac{\dst H_{(3)}}{\dst [H_{(3)}, H_{(3)}]} \ar[u]^{D_{H_{(3)}}}_{\cong} \ar[r] & \rule{0mm}{20pt}\rule[-14pt]{0mm}{20pt}\frac{\dst H_{(3)}H_2}{\dst H_2} \ar@{^(->}[uu] \ar[r] & 0  
  }\EE
 As both rows are exact, we see that 
 \BE\label{intersec}
 \frac{\dst I(G)I(H)}{\dst I^2_{\calG}(G)I(H) + I(G)I^2_{\calH}(H)}\hspace{7pt}\cap \hspace{7pt} \Imm{j_1} = 
  \frac{\dst I(H_{(3)}\cap H_2
  )+I^2_{\calG}(G)I(H) + I(G)I^2_{\calH}(H)}{\dst I^2_{\calG}(G)I(H) + I(G)I^2_{\calH}(H)}
  \EE
 Now consider the following commutative diagram where  $\mu((gG_{(2)})\ot (hH_{(2)})) = (g-1)(h-1)+I^2_{\calG}(G)I(H)$, and similarly for $\mu'$; $D_2(h'H_3)=h'-1+I^2_{\calG}(G)I(H)$ for $h\in H_2$, and similarly or $D_2'$.

\[ \xymatrix{
\rule{0mm}{14pt}\rule[-6pt]{0mm}{0pt}\hspace{14pt} H\ab \sm H\ab \hspace{14pt}\ar[r]^-{(q_2^{\calG}\iota\ab\ot 1)l_2^{\gamma}\hspace{2mm}} \ar@{->>}[d]^-{c_2^{\gamma}} & \rule{0mm}{14pt}\rule[-6pt]{0mm}{0pt}\hspace{14pt}
G\AB \ot H\ab  \ar@{->>}[r]^-{1\ot q_2^{\calH}} \ar@{->>}[d]^-{\mu}
& \rule{0mm}{14pt}\rule[-6pt]{0mm}{0pt}G\AB \ot H\AB \ar@{->>}[d]^-{\mu'}\\
\rule{0mm}{20pt}\rule[-14pt]{0mm}{20pt}   \frac{\dst   H_2}{ \dst H_3} \ar[r]^-{D_2} 
   \ar@{->>}[d]^-{\pi_2}&
\rule{0mm}{20pt}\rule[-14pt]{0mm}{20pt}  \frac{\dst I(G)I(H)}{\dst I^2_{\calG}(G)I(H)} \ar@{->>}[r]^-{\pi_1}  & \rule{0mm}{20pt}\rule[-14pt]{0mm}{20pt}\frac{\dst I(G)I(H)}{\dst I^2_{\calG}(G)I(H) + I(G)I^2_{\calH}(H) }\ar@{->>}[d]^-{\pi_3}\\
  \rule{0mm}{20pt}\rule[-14pt]{0mm}{20pt} \frac{\dst   H_2}{ \dst H_{(3)}} \ar[rr]^-{D_2'}
&  & \rule{0mm}{20pt}\rule[-14pt]{0mm}{20pt}\frac{\dst I(G)I(H)+\Z(G)I(H_{(3)})}{\dst I^2_{\calG}(G)I(H) + I(G)I^2_{\calH}(H)+\Z(G)I(H_{(3)})}
}\]

The upper left-hand square is a pushout by \cite[Theorem 3.6]{Fox2}; the  upper right-hand square is a pushout since $\pi_1$ is surjective and $\Ker{\pi_1} = \mu\,\Ker{1\ot q_2^{\calH}}$. By gluing of pushouts it follows that the upper exterior rectange  also is a pushout. Furthermore, the lower rectangle is a pushout since $\pi_3$ is surjective and $\Ker{\pi_3} = D_2\pi_1\Ker{\pi_2}$ by \REF{intersec}, so again by gluing the two pushout rectangles the exterior square of the whole diagram also is a pushout. Thus
\[ \Ker{\pi_3\mu'} = (q_2^{\calG}\iota\ab\ot q_2^{\calH})l_2^{\gamma} \,\Ker{\pi_2c_2^{\gamma}}\]
We can simplify this description of $\Ker{\pi_3\mu'}$ by observing that the following diagram commutes:
\[ \xymatrix{
H\ab\sm H\ab \ar[rr]^-{(q_2^{\calG}\iota\ab\ot q_2^{\calH})l_2^{\gamma}} \ar[d]^-{\pi_2c_2^{\gamma}} \ar@{->>}[rrd]^-{q_2^{\calH} \sm q_2^{\calH}} & & G\AB \ot H\AB\\
H_2/H_{(3)} \ar@{^(->}[r] & H_{(2)}/H_{(3)}& H\AB \sm H\AB \ar[l]^-{c_2^{\calH}} \ar[u]_-{(\iota\AB\ot 1)l_2^{\calH}}
}\]
 It follows that 
 \BE\label{Kerpi3must} \Ker{\pi_3\mu'} = (\iota\AB\ot 1)l_2^{\calH} \Ker{c_2^{\calH}}
 \EE
Consider the following diagram of plain arrows.

 \N\makebox[14.7cm]{ \makebox[0mm]{
\begin{minipage}{17cm}\small

\[\xymatrix{
\rule{0mm}{20pt}\rule[-14pt]{0mm}{20pt}\hspace{6pt}\frac{\dst  I(G)I(H) + \Z(G)I(H_{(3)})}{\dst {\cal F}_3} \hspace{6pt}\ar@{^(->}[r]^-{j_2} \ar@{=}[d]& \hspace{6pt}
\rule{0mm}{20pt}\rule[-14pt]{0mm}{20pt}
\frac{\dst {\cal F}_2}{\dst {\cal F}_3} \hspace{6pt}\ar@{->>}[r]^-{\pi_4} \ar@{^(->}[d]^-{j_3}& \hspace{6pt}
\rule{0mm}{20pt}\rule[-14pt]{0mm}{20pt}
\coker{j_2} \hspace{6pt}\ar@{.>}[r]^-{\beta} \ar[d]^-{\overline{\pi_5j_3}} & 
\rule{0mm}{20pt}\rule[-14pt]{0mm}{20pt}
\frac{\dst H_{(2)}}{\dst H_2H_{(3)}}\ar@{^(->}[d]\\
\rule{0mm}{20pt}\rule[-14pt]{0mm}{20pt}
\hspace{6pt}\frac{\dst  I(G)I(H) + \Z(G)I(H_{(3)})}{\dst {\cal F}_3} \hspace{6pt}\ar@{^(->}[r] &
\rule{0mm}{20pt}\rule[-14pt]{0mm}{20pt}
 \hspace{6pt} \frac{\dst \Z(G)I(H)}{\dst {\cal F}_3} \hspace{6pt}\ar@{->>}[r]^-{\pi_5} &  
\rule{0mm}{20pt}\rule[-14pt]{0mm}{20pt}
\frac{\dst \Z(G)I(H)}{\dst I(G)I(H) + \Z(G)I(H_{(3)})} \ar[r]_-{\cong}^-{(\overline{D_H})^{-1}} & 
\rule{0mm}{20pt}\rule[-14pt]{0mm}{20pt}
\frac{\dst H}{\dst H_2H_{(3)}}\\
}\]

\end{minipage}\ruled
}\rule[-11mm]{0mm}{3mm} }\vspace{6mm}
 
As both rows are exact in the second group, the map $\overline{\pi_5j_3}$ induced by ${\pi_5j_3}$ is injective. Moreover, ${\rm Im}\big((\overline{D_H})^{-1} \rond\overline{\pi_5j_3}\,\big) = 
{  H_{(2)}}/{ H_2H_{(3)}}$, whence the map $\beta$ exists and is an isomorphism. Together with \REF{Kerpi3must} we 
 obtain an exact sequence
 \BE\label{F2/F3sequ}
 \Ker{c_2^{\calH}} \xrightarrow{(\iota\AB\ot 1)l_2^{\calH}} G\AB \ot H\AB \xrightarrow{j_2\pi_3\mu'}
 {\cal F}_2/{\cal F}_3 \xrightarrow{\beta\pi_4} H_{(2)}/H_2H_{(3)} \lra 0
 \EE
 As to ${\rm U}_2^{\cal GH}(G,H)$, consider the following diagram 
 \BE\label{U2GHpush} \xymatrix{
H\AB\sm H\AB \ar[r]^-{l_2^{\calH}} \ar[d]^-{c_2^{\calH}} & 
\hspace{7pt} H\AB\ot H\AB \hspace{7pt} \ar[r]^-{\iota^{AB}\ot 1} \ar[d]^-{\mu_2^{\calH}\nu_{11}^{\cal HH}} & 
\hspace{7pt} G\AB\ot H\AB \hspace{7pt}\ar[d]^-{\nu_{11}^{\cal GH}} \\
  H_{(2)}/H_{(3)} \ar[r]^-{\nu_2} & {\rm U}_2{\rm L}^{\calH}(H) \ar[r]^-{u_2} & {\rm U}_2^{\cal GH}(G,H)
  }\EE
 One easily checks that both squares are pushouts, by construction of the enveloping algebra and of 
${\rm U}_2^{\cal GH}(G,H)$, resp., and taking the gradings into account. Thus the exterior rectangle is a pushout, too, whence $\Ker{\nu_{11}^{\cal GH}}=(\iota^{AB}\ot 1)l_2^{\calH}\Ker{c_2^{\calH}}$, and there is an isomorphism $\overline{u_2\nu_2} \colon 
 H_{(2)}/H_2H_{(3)} = \coker{c_2^{\calH}} \mr{\cong} \coker{
\nu_{11}^{\cal GH}}$ induced by $u_2\nu_2$.  These facts can be rewritten as an exact sequence
\BE\label{U2GHsequ}
 \Ker{c_2^{\calH}} \xrightarrow{(\iota\AB\ot 1)l_2^{\calH}} G\AB \ot H\AB \xrightarrow{\nu_{11}^{\cal GH} }
{\rm U}_2^{\cal GH}(G,H) \xrightarrow{(\overline{u_2\nu_2})^{-1}\pi_6} H_{(2)}/H_2H_{(3)} \lra 0
 \EE
where $\pi_6\colon {\rm U}_2^{\cal GH}(G,H)\auf \coker{
\nu_{11}^{\cal GH}}$ is the quotient map. Now comparing sequences \REF{F2/F3sequ} and \REF{U2GHsequ} we see that it suffices to show
that the map $\theta^{\cal GH}_2\colon {\rm U}_2^{\cal GH}(G,H)$ $ \lra {\cal F}_2/{\cal F}_3$ commutes with the incoming and outgoing maps, since then the five-lemma implies that it is an isomorphism. But the identity $\theta^{\cal GH}_2 \nu_{11}^{\cal GH} = j_2\pi_3\mu'$ is readily checked by going through the definitions; and to check the identity
\BE\label{outgoing}
\beta \pi_4 \theta^{\cal GH}_2  = (\overline{u_2\nu_2})^{-1} \pi_6
\EE
first note that  the two pushout squares in \REF{U2GHpush} imply that
\begin{eqnarray*}
{\rm U}_2^{\cal GH}(G,H) &=& \Imm{u_2}+ \Imm{\nu_{11}^{\cal GH}} \\
&=& \Imm{u_2\nu_2}+ \Imm{u_2\mu_2^{\cal H}\nu_{11}^{\cal HH}} + \Imm{\nu_{11}^{\cal GH}} \\
&=& \Imm{u_2\nu_2}+  \Imm{\nu_{11}^{\cal GH}}
\end{eqnarray*}
Thus it suffices to check \REF{outgoing} after  precomposition with $u_2\nu_2$ and $\nu_{11}^{\cal GH}$, which is immediate.\hfbox\medskip


As to the map $\theta^{\cal G H}_n$ for $n\ge 3$,
part  of its kernel  is determined in \cite{Fox2} for $\cal H = \gamma$; indeed, all  arguments there remain valid for
arbitray $\cal H$ as above, thus providing an explicitly defined subgroup ${\cal R}_n^{\cal GH}$ of $\UGH{}$ contained in \Ker{
\theta^{\cal GH}}. In particular, one has ${\cal R}_n^{\cal GH}=0$ for $n=1,2$, and ${\cal R}_3^{\cal GH}$ is generated by the
elements
  \BE\label{R3GH}  1 \ot (cH_{(4)})  -  \sum_{q=1}^p  (a_qG_{(3)}) \ot (b_qH_{(2)}) - (b_qG_{(3)}) \ot (a_qH_{(2)})  \EE
where $p\ge 1$, $a_q,b_q \in H \cap G_{(2)}$ such that $c = \prod_{q=1}^p [a_q,b_q] \in H_{(3)}$.
It is shown in \cite{Fox2} that ${\cal R}_n^{\cal GH} \, \bar{\rm U}\LH =0$, so the quotient group
  \[   {\rm \bar{U}}^{\cal GH}(G,H) \stackrel{def}{=}  \UGH{}  \Big/ \ULG{}\,\sum_{n\ge 3} {\cal R}_n^{\calG H}  \]
is a graded $\ULG{}$--${\rm UL}^{\cal H}(H)$-bimodule, and  
  $\theta^{\cal G H}$ induces a surjective homomorphism of graded $\ULG{}$--$\ULH{}$-bimodules
   \[  \sepi{\bar{\theta}^{\calG H}}{ {\rm \bar{U}}^{\calG H}(G,H)\:}{\:{\rm
Gr}^{\cal GH}(\Z(G)I(H))}. \ruled\] 
In particular, $ {\rm \bar{U}}^{\cal GG}(G,G) = {\rm U}^{\cal GG}(G,G)$ and  $\bar{\theta}^{\cal GG} = {\theta}^{\cal GG} =
\theta^{\calG}\mu^{\calG}$, by construction of ${\cal R}_n^{\cal GG}$. While we have no information about $\Ker{\bar{\theta}^{\cal GH}_n}$ for $n\ge 4$ it was computed for
$n=3$ in the cases needed in this paper, namely for $\calH=\gamma$  in \cite{Fox2} and for 
$(H,\calH) =(G,\calG)$  in \cite{Q3}. Both
results involve the   \textit{torsion operator} $\delta_1^{\cal GH}$ described below which basically is the difference between a left and a right connecting homomorphism; it
also describes non trivial torsion relations in the homology of nilpotent groups
(see \cite{GoG}) and thus plays an important role in dimension subgroups \cite{HMP}, whence seems to be quite a fundamental phenomenon.

In fact, all torsion operators in this paper arise from decreasing filtrations 
 $\Delta\,\colon\: \IZ(K) = \Delta_1
\supset \Delta_2 \supset \cdots$ of $I(K)$ for some group $K$ by subgroups $\Delta_i$; such a filtration gives rise to short exact sequences 
\BE\label{Delp} 
0 \hspace{6pt}\to\hspace{6pt} \frac{\Delta_{p+1 }}{\Delta_{p+2}} \hspace{6pt} \stackrel{\alpha}{\hra} \hspace{6pt} \frac{\Delta_{p}}{\Delta_{p+2}} \hspace{6pt}\stackrel{\rho}{\to}\hspace{6pt} \frac{\Delta_{p}}{\Delta_{p+1}} \hspace{6pt}\to\hspace{6pt} 0
\EE  
for $p\ge 1$ where $\alpha,\rho$ denote the canonical injection and quotient map, resp.

Moreover, we need Passi's \textit{polynomial groups with respect to ${\cal G}$} which are denoted by 
\[ P_n^{\calG}(G) = I_{\cal G}(G)/I_{\cal G}^{n+1}(G) \]
see \cite{Pa}. 
By Theorem \ref{U1-U2} sequence \REF{Delp} for $(K,\Delta_i) = (G,I^i_{\calG}(G))$  and $p=1$ gives rise to a natural exact sequence
\BE\label{U2P2G} 0 \to {\rm U}_2\LG{} \hspace{1mm} \mr{\bar{\mu}_2^{\calG}} \hspace{1mm}  P_2^{\calG}(G)  \hspace{1mm} \mr{\rho^{\calG}_2}
\hspace{1mm}  G^{AB} \to 0 \:.\EE 
with $\bar{\mu}_2^{\calG}= \alpha \theta_2^{\calG}$ and $\rho^{\calG}_2 = (\theta_1^{\calG})^{-1}\rho$.
Tensoring this sequence by $H\AB$ and the analogous sequence for $H$ by $G\AB$  gives rise to natural exact sequences
  \BE\label{deftauG}
 {\rm Tor}_1^{\mathbb{Z}} (G^{AB},H\AB)  \hspace{1mm} \mr{\tau_{\calG}}  \hspace{1mm} {\rm U}_2\LG{} \ot \ULH{1}  \hspace{1mm}
\mr{\bar{\mu}_2^{\calG} \ot id} 
\hspace{1mm} P_2^{\calG}(G) \ot H^{AB}
 \hspace{1mm} \mr{\rho^{\calG}_2 \ot id}  \hspace{1mm} G^{AB} \ot H\AB \to 0 \ruled \EE
\BE\label{deftauH} 
{\rm Tor}_1^{\mathbb{Z}} (G^{AB},H\AB)  \hspace{1mm} \mr{\tau_{\calH}}  \hspace{1mm} {\rm U}_1\LG{} \ot \ULH{2}  \hspace{1mm}
\mr{id \ot \bar{\mu}_2^{\calH}} 
\hspace{1mm} G^{AB} \ot  P_2^{\cal H}(H)  
 \hspace{1mm} \mr{id \ot \rho^{\calH}_2}  \hspace{1mm} G^{AB} \ot H^{AB} \to 0 \ruleu\ruled \EE
see \cite[Theorem V.6.1]{ML}. Then define the torsion operator
  \BE\label{del1def} \delta_1^{{\cal GH}} = \nu_{12} \tau_{\calH} - \nu_{21} \tau_{\calG} \hspace{2mm}\colon\hspace{2mm} {\rm Tor}_1^{\mathbb{Z}}
(G^{AB},H\AB)
\hspace{2mm}\lra
\hspace{2mm} {\rm U}_{3}^{{\cal GH}}(G,H) \:.\EE

\N To describe $\delta_1^{{\cal GH}}$ more explicitly we recall from \cite[V.6]{ML} the description of explicit canonical generators of the torsion
product $\Tor{A}{B}$ of abelian groups $A,B$. Suppose that $A=A_1/A_2$ and $B=B_1/B_1$ with $A_1,B_1$ some (non necessarily abelian)
groups. Then these generators are of the form 
\BE\label{torsionML} \langle aA_2,k,bB_2 \rangle \hspace{2mm} \mbox{with \hspace{0.0mm} $ a\in A_1$, $k\in \Z$, $b\in B_1$ such that
$a^k\in A_2$ and $b^k\in B_2$.} \EE 
As a model for all subsequent calculations, we give a  detailed computation of $\delta_1^{{\cal GH}}$, as follows.

\begin{lem} \label{del1comp}
Let $\langle  \bar{g},k,\bar{h} \rangle $ be a canonical generator of ${\rm Tor}_1^{\mathbb{Z}}
(G^{AB},H\AB)$. Then
   \BE\label{del1expl} \delta_1^{{\cal GH}} \langle  \bar{g},k,\bar{h} \rangle   =    \bar{g} \ot (h^kH_{(3)})  -  (g^kG_{(3)}) \ot \bar{h}    
+ {k
\choose 2} ( \bar{g}^2 \ot \bar{h} -   \bar{g} \ot \bar{h}^2) \EE
 \N where we recall that the squares $\bar{g}^2,\bar{h}^2$ are taken in the rings $\ULG{}$ and $\ULH{}$, resp.
\end{lem}

\proof We give the computation of $\nu_{21}\tau_{\calG}$, the one of $\nu_{12}\tau_{\calH}$ being symmetric to this one. 

Throughout in this paper, we use the following fundamental relations in $\Z(G)$, for $a\in G_{(i)}$, $b\in G_{(j)}$, $k\in \Z$ and $n\ge 2$ where $[a-1,b-1]$ denotes the ring commutator in $\Z(G)$.
\begin{eqnarray}
ab-1 &=& (a-1) + (b-1) + (a-1)(b-1)\nonumber\\
&\equiv& (a-1) + (b-1) \hspace{6pt}\bmod\hspace{4pt} I_{\calG}^{i+j}(G)
\label{ab-1}
\end{eqnarray}
\begin{eqnarray} 
 a^k -1 &=& (1+(a-1))^k -1 \nonumber\\
&\equiv& \sum_{p=1}^{n-1} {k \choose p} (a-1)^p \hspace{6pt}\bmod\hspace{4pt} I^{ni}_{\calG}(G) \label{binom}
\end{eqnarray}
\begin{eqnarray} 
 [a,b] -1 &=& [a-1,b-1]a^{-1}b^{-1} \nonumber\\
 &=&  [a-1,b-1] + [a-1,b-1](a^{-1}b^{-1}-1)\nonumber\\
&\equiv&  [a-1,b-1] \hspace{6pt}\bmod\hspace{4pt} I^{i+j+1}_{\calG}(G) \label{comm}
\end{eqnarray}
Now using again \cite[Theorem V.6.1]{ML} one has
 \begin{eqnarray*}
\nu_{21}\tau_{\calG}\langle  \bar{g},k,\bar{h} \rangle &=&
\nu_{21} \left( (\alpha \theta^2_{\calG})^{-1} 
\left(k.\rho^{-1}\theta^1_{\calG}(\bar{g})\right)  \ot \bar{h}\right)\\
&=& \nu_{21} \left( (\theta^2_{\calG} )^{-1}
\alpha^{-1}\Big(k.\rho^{-1}\left( (g-1)+I_{\calG}^2(G) \right)
\Big)  \ot \bar{h}\right)\\
&=& \nu_{21} \left( 
(\theta^2_{\calG} )^{-1}\alpha^{-1}\Big(k(g-1)+I_{\calG}^3(G)
\Big)  \ot \bar{h}\right)\\
&=& \nu_{21} \Big( 
(\theta^2_{\calG})^{-1}\Big((g^k-1)-{k\choose 2}(g-1)^2+I_{\calG}^3(G)
\Big) \ot \bar{h}\Big)\quad\mbox{by \REF{binom}}\\
&=& \nu_{21} \Big( 
\Big((g^kG_{(3)})-{k\choose 2}(gG_{(2)})^2
\Big) \ot \bar{h}\Big)\\
&=&  
(g^kG_{(3)})\ot \bar{h} -{k\choose 2}(gG_{(2)})^2
\ot \bar{h} 
\end{eqnarray*}
as desired.\hfbox\bigskip

 If $(H,{\cal H}) = (G, {\cal G})$ we write $\delta_1^{{\cal G}}=\mu_3^{\calG}\delta_1^{{\cal GG}}$.
This map
completely determines  the structure of $Q_3^{\calG}(G)$, see \cite{Q3}.
\vspace{2mm}

\begin{satz}\label{Q3} The following natural sequence of homomorphisms is exact:
  \[ {\rm Tor}_1^{\mathbb{Z}}(G^{AB},G\AB) \hmr{\delta_1^{{\cal G}}}  \ULG{3}  \hmr{\theta_3^{\calG}} Q_3^{\calG}(G) \hmr{} 0 \:.\]
\end{satz}

For $H\neq G$, or $H=G$ but $\calH\neq \calG$, however, the structure of \Ker{\theta_3^{\calG \calH}} is more complicated, see  
\cite{Fox2} for the case $\calH = \gamma$; in the special case of interest here it is described in the next section. The structure of the related groups
${\rm U}_n^{{\cal N}\gamma}(N,N)$ for $n=2,3$ is determined by the following result, cf.\
\cite[Proposition 5.2]{Fox2}.  Recall that here $N\AB= N/[N,G] = N/[N,N][N,T]$, and let
$l_2^{{\cal N}\gamma} = (q_2^{\cal N}\ot id)l_2^{\gamma}\,\colon\, N\ab \sm N\ab \to N\AB\ot N\ab $. 
\vspace{2mm}

\begin{prop}\label{UnGH} There are canonical isomorphisms
  \[ {\rm U}_2^{{\cal N}\gamma}(N,N) \hspace{2mm}\cong\hspace{2mm}   N\AB\ot N\ab \Big/ l_2^{{\cal N}\gamma}\Ker{c_2^{\gamma}} 
\]
\[ \ruleu {\rm U}_3^{{\cal N}\gamma}(N,N) \hspace{2mm}\cong\hspace{2mm} {\rm coker}\Big( \epsilon = \left( 
\begin{matrix} 
  \phantom{-} c_2^{\cal N} \ot id &   0  &    0  \cr
  - l_2^{\cal N} \ot id  &   i^{{\cal N} \gamma\gamma}   &    i^{\gamma\gamma\gamma} 
\end{matrix}\right) \Big)\:, \]

  \[\begin{matrix} 
\ruled
( (N^{AB} \sm N^{AB}) \ot N\ab ) \hspace{4mm}\oplus\hspace{4mm} ( N^{AB} \ot  l_2^{\gamma} \Ker{c_2^{\gamma}} )
\hspace{4mm}\oplus\hspace{4mm}
l_{33}^N \Ker{c_{33}^N} \cr
\mapdown{\epsilon }  \cr
\rule{0mm}{6mm}\ruled
(  N_{(2)}/ N_{(3)} \ot N\ab )  \hspace{4mm}\oplus\hspace{4mm}  ( N^{AB} \ot N^{AB} \ot N\ab )
\end{matrix}\]

\N whose inverse maps are induced by multiplication in $\Z(G)$. Here
 the homomorphisms ${\rm L}_3^{\gamma}(N)  \ml{c_{33}^N}  (N\ab)^{\ot 3}  \mr{l_{33}^N}  (N\ab)^{\ot 3}$ are defined such
that for
$x,y,z\in N\ab$, 
$c_{33}^N(x\ot y
\ot z)$ is the triple Lie bracket $[x,[y,z]]$ in the Lie algebra ${\rm L}^{\gamma}(N)$ and  $l_{33}^N(x\ot y \ot
z)$ is the triple Lie bracket $[x,[y,z]]$ in the tensor algebra $T(N\ab)$. Furthermore, we note $i^{{\cal N} \gamma\gamma} =  id \ot
q_2^{\cal N}
\ot id\,\colon\,N^{AB} \ot N\ab \ot N\ab \to N^{AB} \ot N^{AB} \ot N\ab $ 
 and
$i^{\gamma\gamma\gamma} = q_2^{\cal N} \ot q_2^{\cal N} \ot id\,\colon\, N\ab \ot N\ab \ot N\ab  \to  N^{AB} \ot N^{AB} \ot N\ab $.
\end{prop}\vspace{3mm}

Actually, the first isomorphism can be easily deduced from the fact that the exterior rectangle in diagram \REF{U2GHpush} is a pushout: for $\calH=\gamma$, $c_2^{\calH}$ is surjective, hence so is $\nu_{11}^{\cal GH}$. Consequently, the map
\BE\label{U2Ggamma}
\overline{\nu_{11}^{\cal G\gamma}} \,\colon\, \frac{G\AB\ot H\ab}{(\iota\AB\ot 1)l_2^{\gamma}\,\Ker{c_2^{\gamma}}} \hspace{6pt} \mr{\cong} \hspace{6pt}{\rm U}_2^{{\cal G\gamma} }(G,H)
\EE
induced by $\nu_{11}^{\cal G\gamma}$ is an isomorphism.

Moreover, the group ${\rm U}_2^{{\cal N}\gamma}(N,N)$ can be embedded into a natural exact sequence, as follows.
Consider the following part of a 6-term-exact sequence
  \[ {\rm Tor}_1^{\mathbb{Z}}(N\AB,N\AB) \hmr{\tau}  N\AB \ot (N_{(2)}/N_2)  \hmr{id\ot i} N\AB \ot N\ab \hmr{id\ot q_2^{\cal N}} N\AB\ot
N\AB
\to 0\]
and let the map $[\,,\,] \,\colon\, N\AB \ot (N_{(2)}/N_2) \:\to\: N_2 \big/ [N_{(2)}\,,N_{(2)}]N_3$ be induced by the commutator pairing of $N$, so that
$[\,,\,] \tau \langle n_1N_{(2)},k,n_2N_{(2)} \rangle  =  [n_1,n_2^k]\, [N_{(2)}\,,N_{(2)}]N_3$.  Furthermore, it follows from 
commutativity of the exterior rectangle in \REF{U2GHpush} (for $H=G=N$, $\calG =\cal N$ and $\calH =\gamma$) that there is a  homomorphism
\[ \overline{l_2^{\gamma}} \mapco  N_2 \big/ [N_{(2)}\,,N_{(2)}]N_3 \to {\rm U}_2^{{\cal N}\gamma}(N,N)  \]
such that for $n_1,n_2\in N$ one has $\overline{l_2^{\gamma}} ( \overline{[n_1,n_2]} ) = {\nu_{11}^{\cal N,\gamma}} \big((n_1N_{(2)}) \ot (n_2N_{2} ) - (n_2N_{(2)}) \ot (n_1N_{2} )\big)$. Finally, for an abelian group   $A$ and $m\ge 1$ let ${\rm SP}^m(A)= A\htt{m}/\Sigma_m$ denote
the symmetric $m$-fold tensor product, and let $\sigma_m=\sigma_m^A \mapco A\htt{m} \auf {\rm SP}^m(A)$ be the canonical projection.

\begin{prop}\label{U2NNsequ}\quad The following
sequence of natural homomorphisms is exact:
  \BE\label{U2inexsequ} {\rm Tor}_1^{\mathbb{Z}}(N\AB,N\AB) \mr{[\,,\,] \tau}   N_2 \big/ [N_{(2)}\,,N_{(2)}]N_3   \mr{\overline{l_2^{\gamma}}} {\rm U}_2^{{\cal N}\gamma}(N,N) 
 \mr{\overline{\sigma_2(1\ot q_2^{\cal N})}} {\rm SP}^2(N\AB) \to 0\:.\EE
 Moreover, $\overline{\sigma_2(1\ot q_2^{\cal N})}$ has a homomorphic splitting if $q_2^{\cal N}$ does and  $N\AB$ is either finitely generated or uniquely $2$-divisible.
\end{prop}\vspace{1mm}

Note that in the case $\cal N=\gamma$ the assertion reduces to the well-known exact sequence 
\BE\label{U2SP2} 0 \lra N_{(2)}/N_{(3)} \lra {\rm U}_2{\rm L}^{\cal N}
(N) \lra {\rm SP}^2(N\AB) \lra 0\EE
and its splitting property if $N\AB$ is either finitely generated or uniquely $2$-divisible, cf.\ \cite{Pa}. For arbitary $\cal N$, the latter facts can in fact be deduced from the left-hand pushout square in \REF{U2GHpush}, in  the same way as we now deduce the proposition from the exterior pushout  rectangle in  \REF{U2GHpush}.\medskip

\proofof{Proposition \ref{U2NNsequ}}
Recall that the exterior   rectangle in  \REF{U2GHpush} is a pushout. This implies the  identity
$\Ker{u_2\nu_2} = c_2^{\gamma}\,\Ker{(\iota^{AB} \ot 1)l_2^{\gamma}}$ and an isomorphism $\coker
{u_2\nu_2} \hcong \coker{(\iota^{AB} \ot 1)l_2^{\gamma}}$. But here $\iota^{AB}=q_2^{\cal N}$, whence the exactness of sequence \REF{U2inexsequ} in $N_2 \big/ [N_{(2)}\,,N_{(2)}]N_3$ follows from Lemma 2.7 in \cite{D3F2}. Moreover, by right-exactness of the tensor product, 
the kernel of the composite map $N\AB\ot N\ab \Sur{1\ot q_2^{\cal N}} N\AB\ot N\AB \Sur{\sigma_2} {\rm SP}^2(N\AB)$ equals $\Imm{( q_2^{\cal N}\ot 1)l_2^{\gamma}}$, 
whence there is an isomorphism $\coker{(q_2^{\cal N} \ot 1)l_2^{\gamma}} \hcong {\rm SP}^2(N\AB)$ induced by $\sigma_2(1\ot q_2^{\cal N})$.
This implies exactness of sequence \REF{U2inexsequ} in  
${\rm U}_2^{{\cal N}\gamma}(N,N) $. Finally, the splitting assertion follows from the well-known fact that $\sigma_2^A$ splits if $A$ is either finitely generated or uniquely $2$-divisible: if $s_1,s_2$ are splittings of $q_2^{\cal N}$ and of $\sigma_2^{N\AB}$, resp., $\nu_{11}^{\cal N\gamma} (1\ot s_1)s_2$ is a splitting of  $\overline{\sigma_2(1\ot q_2^{\cal N})}$
 .\hfbox

\section{The first four Fox and augmentation quotients
}

Throughout this section $G$ is supposed to be the semidirect product of a normal subgroup $N$ with some subgroup $T$.
The following groups are given a complete functorial description in terms of the Lie algebras $\LN$ and ${\rm L}^{\gamma}(T)$,  for all   $G=N\rtimes T$:
\begin{itemize}
\item the quotients $Q_n(G)$, $Q_n(G,T)$ and $Q_n(G,N)$ for $n\le 3$;

\item   the direct factors $\calK_4 / \calK_5$ of $Q_4(G)$ and $Q_4(G,T)$
and the direct factor $\Gamma_3^*I(N)/\Gamma_4^*I(N)$ of $Q_4(G,N)$, see \REF{AugTallg}, \REF{FoxTallg} and \REF{FoxNallg}.
\end{itemize}

The groups $Q_4^{\calG}(G)$, and hence the direct factors $Q_4^{\cal N}(N)$ and $Q_4^{\gamma}(T)$ of $Q_4(G)$ and $Q_4(G,T)$, 
were determined for finite $G$ by Tahara \cite{Ta4}; so the only term  for $n=4$ whose structure remains almost completely unknown is the direct factor $\Lambda_3I(N)/\Lambda_4I(N)$ of $Q_4(G,N)$ (it is only computed under very restrictive assumptions in Corollary \ref{FoxNUU} below).

The proofs of all results of this section are deferred to section 4.\V

   \comment{

in terms of tensor and torsion products among groups   $\UTn{i}$ and
$\UNn{j}$, or explicit quotients of such terms. The only group not computed in this paper is the direct factor
$\Lambda_3I(N)/\Lambda_4I(N)$ of $Q_4(G,N)$. All proofs are deferred to section 4.\V
}

The groups  $Q_2(G,K)$ were determined by Tahara \cite{Ta} for $K=G$ and by Karan and Vermani \cite{KV29}, \cite{KV30} for
$K=N,T$, after partial results of Khambadkone \cite{Kh13}, \cite{Kh14}; we quote the results here (expressed in the language of enveloping rings) for completeness but also because they are easily reproved using our general approach, see
section 4.
\vspace{2mm}

\begin{satz}\label{Fox2} There are natural isomorphisms
 \[  Q_2(G)    \hspace{2mm}\cong\hspace{2mm}  {\rm U}_2{\rm L}^{\cal N}(N) \hspace{2mm}\oplus\hspace{2mm} 
{\rm U}_2{\rm L}(T)  \hoplus  N\AB \ot T\ab  \]
 \[  Q_2(G,T)    \hspace{2mm}\cong\hspace{2mm}    {\rm U}_2{\rm L}(T)  \hoplus  N\AB \ot T\ab  \]
 \[  Q_2(G,N)    \hspace{2mm}\cong\hspace{2mm}   {\rm U}_2^{{\cal N}\gamma}(N,N)  \hoplus  T\ab \ot N\ab
\]
\end{satz}


Our description of $Q_n(G,H)$ for $n=3,4$ below involves various torsion operators  coming from connecting homomorphisms as in \cite[Theorem V.6.1]{ML}. To keep notations simple we denote by $\tau_k^{\ssst \Box}, \hat{\tau}_k^{\ssst \Box}$, $k=1,2$ and $\sst\Box$ some (or no) superscript, a connecting homomorphism induced by a short exact sequence of abelian groups in the $k$-th variable. In particular, for $p,q=1,2$ we have maps \vspace{11pt}

\N\makebox[14.7cm]{ \makebox[0mm]{
\begin{minipage}{17cm}\small
\BE\label{tau1tau2} \frac{\UNn{p+1}}{\Imm{\delta_1^p}} \ot \UTn{q}    \hml{\tau_1^{pq}}   \Tor{\UNn{p}}{\UTn{q}}  \hmr{\tau_2^{pq}}  \UNn{p} \ot \frac{\UTn{q+1}}{\Imm{\delta_2^q}}  \EE
\end{minipage}\ruled
}
}
\vspace{1pt}

with $\delta_1^1,\delta_2^1=0$, $\delta_1^2=\delta_1^{\cal N}$ and $\delta_2^2=\delta_1^{\gamma}$, see \REF{del1def}, where

\begin{itemize}
\item $\tau_1^{1q}$ and $\tau_2^{p1}$
are induced by the short exact sequence  \REF{U2P2G} for $(G,{\cal G}) = (N,{\cal N})$ and\ $(T, \gamma)$, resp.; explicitly, we have
\BE\label{tau1q}
\tau_1^{1q}\langle nN_{(2)},k,x\rangle \:=\: \Big((n^kN_{(3)}) - {k\choose 2}(nN_{(2)})^2 \Big)\ot x
\EE
\BE\label{taup1}
\tau_2^{p1}\langle y,k,tT_2\rangle \:=\: y\ot \Big( (t^kT_3) - {k\choose 2}(tT_2)^2 \Big) 
\EE
for suitable $n,t,x,y,k$; cf.\ the calculation of $\tau_{\calG}$ in the proof of Lemma \ref{del1comp}.

\item $\tau_1^{2q}$ and $\tau_2^{p2}$
are induced by the short exact sequence
  \BE\label{U3U2sequ} 0 \hspace{1mm}\to\hspace{1mm}  \coker{\delta_1^{\cal G}}   \hmr{\bar{\mu}^{{\cal G}}_3} \hspace{2mm}  I^2_{\cal G}(G)/ I^4_{\cal G}(G)  \hmr{\rho_3^{\cal G}} \UGn{2}  \hspace{1mm}\to\hspace{1mm} 0 \EE  
obtained from sequence \REF{Delp} for $p=2$ combined with Theorems  \ref{Q3} and \ref{U1-U2}, for $(G,{\cal G}) = (N,{\cal N})$ and\ $(T, \gamma)$, resp.
\end{itemize}

The maps $\tau_1^{2q}$ and $\tau_2^{p2}$, occuring in Theorem \ref{Fox4} below, can be made explicit if $N\ab$ and $T\ab$ are finitely generated, by the following formula (and its mirror-symmetric version, permuting $G$ and $A$, which gives the corresponding map $\tau_2$).

\begin{prop}\label{tau2q}
Let $G$ be a group with N-series $\calG$ and $A$ be an abelian group such that $G\AB$ is finitely generated. Let ${\rm Tor}(G\AB) = \bigoplus_{i=1}^m \Z/e_i\Z\langle g_iG_{(2)}\rangle$, $g_i\in G$,  be a cyclic decomposition of the torsion subgroup of $G\AB$, and let $e_{ij}$ be the greatest common divisor of $e_i$ and $e_j$, and $p_{ij},q_{ij}\in\Z$ such that $e_{ij}=e_ip_{ij}+e_jq_{ij}$.  Furthermore, let $q_3\colon \UGn{3} \auf \coker{\delta_1^{\cal G}}$ be the quotient map.
Then the connecting homomorphism 
\[\tau_1\,\colon\, \Tor{\ULG{2}}{A} \lra \coker{\delta_1^{\calG}} \ot A\]
induced by sequence \REF{U3U2sequ} is given as follows. Let $\langle x,k,a\rangle$ be a canonical generator of $\Tor{\ULG{2}}{A}$. According to the split exact sequence \REF{U2SP2}, $x$ can be uniquely written in the form $x=(gG_{(3)}) + \sum_{1\le i\le j\le m} l_{ij}(g_iG_{(i)}) \ot (g_jG_{(j)})$ with $g\in G_{(2)}$, $l_{ij}\in \Z$ such that $g^k\in G_{(3)}$ and $e_{ij}$ divise $kl_{ij}$ for all $1\le i\le j\le m$. Then
\[ \tau_1\langle x,k,a\rangle = \T{x}\ot a \]
with 
\begin{eqnarray*}
\T{x}&=& q_3\left( \bigg( (g^kG_{(4)}) + \sum_{1\le i\le j\le m} \frac{kl_{ij}p_{ij}}{e_{ij}} [g_i^{e_i} G_{(3)}, g_j  G_{(2)}]\right.\bigg)\\
& & {} \oplus \sum_{1\le i\le j\le m} \bigg( \frac{kl_{ij}p_{ij}}{e_{ij}} ( g_j  G_{(2)})(g_i^{e_i} G_{(3)}) + \frac{kl_{ij}q_{ij}}{e_{ij}} ( g_i  G_{(2)})(g_j^{e_j} G_{(3)}) \bigg)\\
& & \left. {} \ominus \sum_{1\le i\le j\le m}   \frac{kl_{ij}}{e_{ij}}\bigg( p_{ij} {e_i\choose 2} (g_iG_{(2)})^2( g_j  G_{(2)}) + q_{ij} {e_j\choose 2} (g_iG_{(2)})( g_j  G_{(2)})^2 \bigg) 
\right)
\end{eqnarray*}
Here the symbols $\oplus,\ominus$ mean $+,-$, resp., but also indicate that the three summands they link together lie in the three different direct components of the decomposition 
\[ \UGn{3} \hspace{9pt}\cong\hspace{9pt} 
G_{(3)}/G_{(4)} \hspace{6pt}\oplus\hspace{6pt} (G_{(1)}/G_{(2)})\ot (G_{(2)}/G_{(3)}) \hspace{6pt}\oplus\hspace{6pt} {\rm SP}^3(G_{(1)}/G_{(2)})\]
\end{prop}

We omit the proof since it consists of  a calculation which is  straightforward along the same lines as the computation of $\delta_1^{\cal GH}$ in the proof of Lemma \ref{del1comp} and of 
$\xi_3$ in the proof of Theorem \ref{Fox4} in section 4.\V

Moreover, we throughout identify ${\rm U}_1{\rm L}^{\cal K}(K)$ with $K\AB$ via the isomorphism $\theta_1^{\cal K}$, see Theorem \ref{U1-U2}, for $(K,{\cal K}) = (N,\cal N)$ or $(T,\gamma)$. Moreover, we identify $\bar{\rm U}_2^{{\cal N} \gamma}(N,N)$ with $\frac{N\AB\ot N\ab}{(q_2^{\cal N}\ot 1)l_2^{\gamma}\,{\rm Ker}(c_2^{\gamma})}$ via the isomorphism $\overline{\nu_{11}^{\cal N\gamma}}$ in \REF{U2Ggamma}.

Also recall that the structure of the group $\bar{\rm U}_3^{{\cal N} \gamma}(N,N)= {\rm U}_3^{{\cal N} \gamma}(N,N)/{\cal R}^{N\gamma}_3$ is explicitly given by Proposition \ref{UnGH} and the generators of ${\cal R}^{N\gamma}_3$ described in \REF{R3GH}, but we will merely use the original definition of ${\rm U}_3^{{\cal N} \gamma}(N,N)$ in \REF{UGHdef}.\V

Now we are ready to describe the structure of the groups $Q_3(G,H)$ for $H=G,N,T$.\medskip

\begin{satz}\label{Fox3} The terms on the right hand side of the decompositions
 \[   Q_3(G)    \hspace{2mm}=\hspace{2mm}  Q_3^{\cal N}(N) \hspace{2mm}\oplus\hspace{2mm} 
Q_3^{\gamma}(T)  \hoplus  \calK_3 / \calK_4  \]
 \[  \ruled Q_3(G,T)    \hspace{2mm}=\hspace{2mm}    Q_3^{\gamma}(T)  \hoplus  \calK_3 / \calK_4  \]
 \[  Q_3(G,N)    \hspace{2mm}=\hspace{2mm}    \Lambda_2I(N)/\Lambda_3I(N)  \hoplus  \Gamma_2^*I(N)/\Gamma_3^*I(N) \]
are determined by Theorem \ref{Q3} and the following natural exact sequences.
  \[ \Tor{N\AB}{T\ab} \hmr{\delta_2} \UNn{2} \ot T\ab \hoplus N\AB \ot \UTn{2} \hmr{\mu_2}  \calK_3 / \calK_4 \to 0 \]
  \[ \Tor{T\ab}{N\ab} \hmr{\delta_3}  T\ab  \ot {\rm U}_2^{{\cal N}\gamma}(N,N)
\hoplus \UTn{2} \ot N\ab
 \hmr{\mu_3}  \Gamma_2^*I(N)/\Gamma_3^*I(N) \to 0 \]
  \[\begin{matrix} 
\ruleu\ruled {\rm Tor}_1^{\mathbb{Z}}(N\AB,N\ab) \hspace{4mm}\oplus\hspace{4mm}  
{\rm Ker}\Big( [\,,\,]\tau\,\colon\, {\rm Tor}_1^{\mathbb{Z}}(N\AB,N\AB) \to N_2 \big/ [N_{(2)}\,,N_{(2)}]N_3 \Big)  \cr
\mapdown{(\delta_4,\delta_5)} \cr
\ruleu\ruled \bar{\rm U}_3^{{\cal N} \gamma}(N,N)  \cr
\surdown{\bar{\theta}_3^{{\cal N} \gamma}}  \cr
\ruleu\ruled  \Lambda_2I(N)/\Lambda_3I(N)
\end{matrix}\]
Here the homomorphisms $\mu_2,\mu_3$ are induced by $\theta_2^{\cal N}, \theta_2^{\gamma}, \theta_2^{{\cal N}\gamma}$ followed by 
multiplication in $\Z(G)$,   
$\delta_2 = (-\tau_1^{11},\tau_2^{11})^t$, $\delta_3$, $\delta_4$ are homomorphisms and $\delta_5$ is an additive relation of undeterminacy
\Imm{\delta_4}, defined as follows. Using the identifications in Proposition \ref{UnGH} one has for suitable $n\in N$ and $t\in T$, see \REF{torsionML}:  \begin{eqnarray*}
  \delta_2 \langle n N_{(2)}\,,k,tT_2 \rangle  &=&  \bigg( {}- (n^k N_{(3)}) \ot (tT_2) + {k\choose 2} (n  N_{(2)})^2
\ot (tT_2) \hspace{2mm}, \\  
& & (n  N_{(2)}) \ot (t^k T_3) -  {k\choose 2} (n  N_{(2)})  \ot
(tT_2)^2 \bigg)  
  \end{eqnarray*}
     
       \[  \delta_3 \langle tT_2 \,,k,n N_2 \rangle  =  \bigg( (tT_2) \ot \nu_{11}^{\cal N\gamma} \Big( 
       \sum_{i=1}^p \Big( (n_iN_{(2)}) \ot (n_i\st N_2) -   (n_i\st N_{(2)}) \ot (n_i N_2) \Big) \]
       \[ \hspace{4.5cm}  - {k\choose 2} (n  N_{(2)}) \ot (n  N_2)
  \Big) \hspace{1mm}, \hspace{1mm} \Big( {k\choose 2} (tT_2)^2
     - (t^kT_3) \Big)  \ot (nN_2) \bigg)  \]

\N where  $p\ge 1$ and $n_i,n_i\st \in N$ such that $n^k = \prod_{i=1}^p [n_i,n_i\st]$.  Furthermore, for suitable $a,b\in N$ and denoting by $\pi\colon  {\rm U}_3^{{\cal N} \gamma}(N,N)\auf  \bar{\rm U}_3^{{\cal N} \gamma}(N,N)$ the quotient map,
  \begin{eqnarray*}
  \delta_4 \langle aN_{(2)} \,,k, b N_2 \rangle  &=& \pi\Big( (aN_{(2)}) \ot (b^kN_3)  -  (a^kN_{(3)}) \ot (bN_2)   \\
& &  + {k
\choose 2} \Big( (aN_{(2)})^2 \ot (bN_2) -   (aN_{(2)}) \ot (bN_2)^2 \Big)  \Big)
  \end{eqnarray*}  
Finally,  for $\sum_{r=1}^s \langle a_r N_{(2)} \,,k_r, b_r  N_{(2)} \rangle \in {\rm Tor}_1^{\mathbb{Z}}(N\AB,N\AB)$ such that
 $\prod_{r=1}^s [a_r,b_r^{k_r}] =   e \prod_{q=1}^p [c_q,d_q]$ with $c_q,d_q\in [N,G]$ and $e \in N_3$,
  \begin{eqnarray*}
  \delta_5 \bigg( \sum_{r=1}^s \langle a_r N_{(2)} \,,k_r, b_r  N_{(2)} \rangle  \bigg)  &=& 
  \sum_{r=1}^s (a_r^{k_r}N_{(3)}) \ot (b_r N_2) -  (b_r^{k_r}N_{(3)}) \ot (a_r N_2) \\ 
& - &   \sum_{r=1}^s {k_r \choose
2} \Big(  (a_rN_{(2)}) \Big((a_rN_{(2)}) - (b_rN_{(2)}) \Big) \ot (b_rN_2) \Big) \\
& - &  \sum_{q=1}^p (c_q N_{(3)}) \ot (d_qN_2) - (d_q
N_{(3)}) \ot (c_qN_2)\\
&-& 1 \ot (eN_4)  + \Imm{\delta_4}\:.
  \end{eqnarray*}

\end{satz}

This result generalizes and extends the computation of $Q_3(G)$ for finite $G$ in \cite{Ta}  and of  $Q_3(G,T)$ and $\Gamma_2^*I(N)/\Gamma_3^*I(N)$ for finite $G$ and nilpotent $T$
in \cite{Kh15}, \cite{Kh12}. It seems, however, that the group $\Lambda_2I(N)/\Lambda_3I(N)$ has not been determined before, not even in special cases.

We now turn to the case $n=4$ where, apart from the direct factors $Q_4^{\cal N}(N)$ and $Q_4^{\gamma}(T)$, 
nothing seems to be known unless $N$ and $T$  satisfy  certain torsion-freeness conditions, see section 3.\V

\nc{\hhoplus}{\hspace{1mm}\oplus\hspace{1mm} }

\begin{satz}\label{Fox4} The direct factor $\calK_4/\calK_5$ of $Q_4(G)$ and $Q_4(G,T)$ (see \REF{AugTallg} and \REF{FoxTallg}) is determined by the following tower of successive natural quotients where $\Ker{\pi_k} = \Imm{\xi_k}$, $k=1,2,3$, and 
\[\xi_1=\left(\begin{matrix} \delta_1^{\cal N}\ot 1 & 0 & 0 \cr 0 & 0 & 1\ot  \delta_1^{\cal \gamma}\end{matrix}\right)^t\:,\quad\xi_2 = \left(\begin{matrix} {} - {\tau}^{21}_1 & \tau_2^{21} & 0 \cr
 0 & {} -\tau_1^{12} &  {\tau}_2^{12}\end{matrix}\right)^t \:,\]
cf.\ the explicit description of these maps in \REF{tau1q}, \REF{taup1} and Proposition \ref{tau2q}.

\N\makebox[14.7cm]{ \makebox[0mm]{
\begin{minipage}{17cm}\small

\[ \begin{matrix}
\left.\begin{matrix} \ruled \Tor{N\AB\hspace{-1mm}}{N\AB} \ot T\ab \cr \oplus\hspace{1mm}
 N\AB \ot \Tor{T\ab}{T\ab} \end{matrix} \right\}
 &  \mr{\xi_1}  &  \UNn{3} \ot T\ab \hoplus \UNn{2} \ot \UTn{2} 
 \hoplus N\AB \ot \UTn{3} \cr
  & & \surdown{\pi_1} \cr
\ruleu\ruled
\left.\begin{matrix} \ruled\Tor{\UNn{2}}{T\ab} \cr  \oplus \hspace{1mm}\Tor{N\AB\hspace{-1mm}}{\UTn{2}}\end{matrix} \right\}
  & \mr{\xi_2}  &  \coker{\delta_1^{\cal N}}  \ot T\ab \hoplus \UNn{2} \ot \UTn{2} 
 \hoplus N\AB \ot \coker{\delta_1^{\gamma}} \cr
   & & \surdown{\pi_2} \cr
\ruleu\ruled \Ker{\delta_2} & \mr{\xi_3}  &  \coker{\xi_2} \cr
      & & \surdown{\pi_3} \cr
  & & \ruleu \calK_4/\calK_5
\end{matrix}\]

\end{minipage}\ruled
}\rule[-11mm]{0mm}{3mm} }\vspace{6mm}
To describe $\xi_3$ we here suppose that $N$ and $T$ are finitely generated  with cyclic decompositions of the torsion subgroups {\rm Tor}$(N\AB) = \bigoplus_{i=1}^r \Z/a_i\Z \langle n_iN_{(2)} \rangle$
and {\rm Tor}$(T\ab) = \bigoplus_{j=1}^s \Z/b_j\Z \langle t_jT_2 \rangle$.
Let $d_{ij}$ be the greatest common divisor of $a_i$ and $b_j$, and let $p_{ij},q_{ij} \in \Z$ such that $d_{ij}=a_ip_{ij} + b_jq_{ij}$.
Then an element $\omega = \sum_{i,j} \langle n_iN_{(2)}, k_{ij}, t_jT_2 \rangle \in\Tor{N\AB}{T\ab} $ lies in $\Ker{\delta_2}$ if and only if the following three conditions (i) - (iii) are satisfied:

\begin{rom}
\item \hspace{2mm} $\forall 1\le i\le r$, $\forall 1\le j\le s$, \hspace{5mm} $\frac{k_{ij}}{d_{ij}}$ is even if $k_{ij}$ is even;

\item \hspace{2mm} $\forall 1\le i\le r$,  \hspace{5mm} $\prod_{j=1}^s t_j^{k_{ij}} = u_i^{a_i} v_i$ \hspace{2mm}with $u_i\in T_2$ and $v_i\in T_3$;

\item  \hspace{2mm} $\forall 1\le j\le s$, \hspace{5mm} $\prod_{i=1}^r n_i^{k_{ij}} = y_j^{b_j} z_j$ \hspace{2mm}with $y_j\in N_{(2)}$ and $z_j\in N_{(3)}$.
  
\end{rom}


 
 \N In this case, 
 \begin{eqnarray*}
  \xi_3(\omega)  &=& \pi_2\pi_1 \left( {} - \sum_{j=1}^s \bigg( (z_j N_{(4)}) - \sum_{i=1}^r \bigg(
   \frac{p_{ij}}{d_{ij}} {k_{ij} \choose 2}   (n_i^{a_i}N_{(3)})(n_iN_{(2)})
    \right. \\  
  & + & \left.  \left.  \left( {k_{ij} \choose 3}  -  \frac{p_{ij}}{d_{ij}} {k_{ij} \choose 2} {a_i \choose 2}\right)   (n_iN_{(2)})^3        \right)\right)
     \ot (t_jT_2) \hspace{1mm}\mathbb{,}\hspace{1mm} \\
 & &    \sum_{i=1}^r \bigg(   (n_i^{a_i}N_{(3)}) - {a_i \choose 2}  (n_iN_{(2)})^2 \bigg)
          \ot   \bigg( (u_iT_3) -  \sum_{j=1}^s \frac{p_{ij}}{d_{ij}} {k_{ij} \choose 2} (t_jT_2)^2 \bigg) \\
& - &   \sum_{j=1}^s  \bigg( (y_jN_{(3)}) - \sum_{i=1}^r   \frac{q_{ij}}{d_{ij}} {k_{ij} \choose 2} (n_iN_{(2)})^2 \bigg) 
       \ot \bigg( (t_j^{b_j}T_3) - {b_j \choose 2} (t_jT_2)^2 \bigg) \hspace{1mm}\mathbb{,}\hspace{1mm} \\
& &   \sum_{i=1}^r (n_iN_{(2)}) \ot  \left( \rule{0mm}{20pt}
\phantom{{ \choose }} \makebox[3mm]{$\phantom{{ \choose }}$} \right.  
\hspace{-30pt}(v_iT_4)
     -   \sum_{j=1}^s \bigg( \frac{q_{ij}}{d_{ij}} {k_{ij} \choose 2}   (t_j^{b_j}T_3)(t_jT_2)  \\
   & + & \left. \left.\left.  \left( {k_{ij} \choose 3} -  \frac{q_{ij}}{d_{ij}} {k_{ij} \choose 2}{b_j \choose 2} \right)(t_jT_2)^3    \right) \rule{0mm}{20pt}\right)\rule{0mm}{23pt}\right)
    \end{eqnarray*}

\end{satz}\vspace{15pt}

Similarly, the direct factor
$\Gamma_3^*I(N)/\Gamma_4^*I(N)$ of $Q_4(G,N)$ can be computed by combining Theorems \ref{KD4/KD5}, \ref{U1-U2}, \ref{Q3} and Remark \ref{symbem}, but the resulting description is considerably more complicated than the one of $\calK_4/\calK_5$ above, so we leave it to the interested reader to write it out.

In principle, one can use the key Proposition \ref{amalgam} to go on and determine $\calK_n/\calK_{n+1}$ and $\Gamma_{n-1}^*I(N)/\Gamma_n^*I(N)$ for   $n\ge 5$, in terms of iterated amalgamations of tensor products of the augmentation quotients of $N$ and $T$ along certain torsion groups, but the results getting more and more complicated we do not attempt to make this explicit. When all these torsion terms vanish, however, the amalgamations degenerate to  neat  direct sum decompositions; this is described in the next section.

  \comment{
Then there is a natural epimorphism
  \[ \theta^{\caln N}_n\,\colon\: \UNN{n} \Sur{} \Lambda_{n-1}I(N)/\Lambda_nI(N) \]
defined as follows: for $n=i+j$ with $i\ge 0$ and $j\ge 1$, $x\in \ULN{i}$, $y\in {\rm U}_j{\rm L}(N)$, $x\st \in \Lambda_i$, $y\st
\in I^j(N)$ such that $\theta^{N,\caln}_i(x) = x\st + \Lambda_{i+1}$ and $\theta^{n,\caln}_j(x) = x\st + I^{j+1}(N)$ one has 
$\theta^{\caln N}_n(x\ot y) = x\st y\st + \Lambda_nI(N)$. A subgroup ${\cal R}_n^{\caln N}$ contained in the kernel of
$\theta^{\caln N}_n$ is explicitly described in \cite{Fox2};  ${\cal R}_3^{\caln N}$ is generated by the elements  
  \[ 1\ot (cN_4) - \sum_{q=1}^p (a_qN_{(3)}) \ot (b_qN_2) - (b_qN_{(3)}) \ot (a_qN_2)  \]
where $p\ge 1$, $a_q,b_q\in N_{(2)}$ such that $c = \prod_{q=1}^p [a_q,b_q] \in N_3$. Let $\bar{{\rm U}}^{\cal N}_{3}(N,N) =
\UNN{3}/{\cal R}_3^{\caln N}$ and $\bar{\theta}^{\caln N}_3\,\colon \bar{{\rm U}}^{\cal N}_{3}(N,N) \auf
\Lambda_{2}I(N)/\Lambda_3I(N)$ be induced by $\theta^{\caln N}_3$. Finally, we abreviate $N\AB = N/N_{(2)}$.\V
}

\section{Fox and augmentation quotients under torsion-freeness assumptions}

Supposing one or more among the groups $N\ab$, $N\AB$ and $T\ab$ torsion-free the groups $Q_n(G,H)$ for $H=G,N,T$ and $n\le 4$ were determined   by Karan and Vermani, see the precise citations below. We here generalize their results to all $n\ge 1$, and improve them by expressing most of the involved groups in terms of enveloping algebras. All proofs are deferred to section 5.\V

We formally put $T_0=T$ and $N_{(0)}=N$.\V

\begin{satz}\label{Aug+FoxT}\quad Let $n\ge 2$. Suppose that there  exists  $k$, $0 \le k \le n-2$,
such that $T_s/T_{s+1} $ is torsion-free for $0\le s \le k$ and that $N_{(t)}/
N_{(t+1)}$ is torsion-free for $0\le t \le n-k-2$. 
Then there are  natural isomorphisms
  \[  Q_n(G)\hspace{2mm}\cong\hspace{2mm}  Q_n^{\cal N}(N) \hspace{2mm}\oplus \hspace{2mm} Q_n(T)
\hspace{2mm}\oplus\hspace{2mm} \bigoplus_{i=1}^{n-1} \hspace{2mm}Q_i^{\cal N}(N) \ot  Q_{n-i}(T) \] 

  \[ Q_n(G,T)\hspace{2mm}\cong\hspace{2mm}  Q_n(T)\hspace{2mm}\oplus\hspace{2mm}
\bigoplus_{i=1}^{n-1} \hspace{2mm}Q_i^{\cal N}(N) \ot  Q_{n-i}(T)   \]
\end{satz}

This implies the results in \cite{KV30}, \cite{KV32} for $n=3$ and in \cite{KV31}, \cite{KV47} for $n=4$ (which correspond to the case $k=1$).\V

\begin{kor}\label{Aug+FoxTU} \quad Suppose that  $T_s/T_{s+1} $ and
$N_{(s)}/ N_{(s+1)}$ are torsion-free for  $1\le s\le n$. Then there are  natural isomorphisms
    \[  Q_n(G) \hspace{2mm}\cong\hspace{2mm}   \bigoplus_{i=0}^{n} \hspace{2mm} 
{\rm U}_i{\rm L}^{\cal N} (N) \ot {\rm U}_{n-i}{\rm L}^{\gamma}(T) \:.\]
     \[ Q_n(G,T) \hspace{2mm}\cong\hspace{2mm}  \bigoplus_{i=0}^{n-1} \hspace{2mm}
{\rm U}_i{\rm L}^{\cal N}(N) \ot {\rm U}_{n-i}{\rm L}^{\gamma}(T) \:,\]
\end{kor}

Using the Poincar\'e-Birkhoff-Witt theorem one deduces from \ref{Aug+FoxTU} the following result which generalizes a theorem of
Sandling and Tahara \cite{Sa-Ta} on the augmentation quotients of an arbitrary group (which can be recovered here by taking
$N=\{1\}$). Recall that  by convention ${\rm SP}^0(X)=\Z$.

\begin{kor}\quad If for  $1\le s\le n$  the abelian groups $T_s/T_{s+1} $ and
$N_{(s)}/ N_{(s+1)}$ are free (in particular if they are torsion-free and $N$ and $T$ are finitely generated) then
\[  Q_n(G)\:\cong\:  \bigoplus_{{\cal I}_1} \bigotimes_{p=1}^{n }
{\rm SP}^{r_p}(N_{(p)}/ N_{(p+1)}) \ot \bigotimes_{q=1}^{n } {\rm SP}^{s_q}(T_q/
T_{q+1}) \:,\]
  \[ Q_n(G,T) \:\cong\:  \bigoplus_{{\cal I}_2} \bigotimes_{p=1}^{n-1} {\rm SP}^{r_p}(N_{(p)}/
N_{(p+1)}) \ot \bigotimes_{q=1}^{n } {\rm SP}^{s_q}(T_q/
T_{q+1}) \:,\]
\ruleu  where the index sets ${{\cal I}_1}$ and ${{\cal I}_2}$ are given by
  \[  {{\cal I}_1} = \bigg\{(r_1, \ldots,r_{n},s_1,\ldots,s_{n}) \,\Big|\,0 \le r_1,
\ldots,
r_{n}\,, s_1,\ldots,s_{n } \le n \mbox{\hspace{1mm} and\hspace{1mm} }
\sum_{p=1}^{n }  r_p p + \sum_{q=1}^{n }  s_q q   = n\, \bigg\} \]
\[ {{\cal I}_2} =  \bigg\{(r_1, \ldots,r_{n-1},s_1,\ldots,s_{n}) \,\Big|\,
0 \le r_1,\ldots,r_{n-1} \le n-1, \hspace{2mm}0 \le s_1,\ldots,s_{n } \le n\, , \]
  \[ \hspace{52mm} \sum_{q=1}^{n}  s_q q \ge 1 \mbox{\hspace{2mm} and\hspace{2mm} }
\sum_{p=1}^{n-1}  r_p p + \sum_{q=1}^{n}  s_q q  = n\, \bigg\}\:.\] 
\end{kor}\vspace{2mm}

  As to the quotients $Q_n(G,N)$ we have the following results.\vspace{2mm}

\begin{satz}\label{FoxN}\quad Suppose that ${\rm
Tor}_1^{\mathbb{Z}}\big(I^i(T)/I^{i+1}(T)\,,\,I(N)/\Lambda_{n-i-1}  I(N) \big) =0$ for  $1\le
i\le n-2$.  Then there is a natural isomorphism
 \[Q_n(G,N)\hspace{2mm}\cong\hspace{2mm} \frac{\dst \Lambda_{n-1}I(N)}{\dst \Lambda_{n}I(N)}
\hspace{2mm}\oplus\hspace{2mm} \bigoplus_{i=1}^{n-1} \hspace{2mm}
Q_i(T) \ot\hspace{1mm}
\frac{\dst \Lambda_{n- i-1}I(N)}{\dst
\Lambda_{n- i}I(N)}  \:.\]
\end{satz}\vspace{2mm}

For $n=3$ this reproduces the main result in \cite{KV30}, \cite{KV32}; for $n=4$ it implies the main result in \cite{KV31}, \cite{KV47}
 since torsion-freeness of $T\ab\cong I(T)/I^2(T)$ and $N\ab\cong I(N)/\Lambda_1I(N)$ imply triviality of our torsion group for $i=1,2$, resp.\vspace{2mm}

\begin{kor}\label{FoxNU}\quad If $T_s/T_{s+1}$ is torsion-free for
$1\le s\le n$ then there is a natural isomorphism
  \[Q_n(G,N) \hspace{2mm}\cong\hspace{2mm}
\bigoplus_{i=0}^{n-1} \hspace{2mm}
{\rm U}_i{\rm L}^{\gamma} (T) \ot\hspace{1mm}
\frac{\dst \Lambda_{n- i-1}I(N)}{\dst
\Lambda_{n- i}I(N)}  \:.\]
\end{kor}

  Now using \cite[Proposition 2.1]{Fox2} we get

\begin{kor}\label{FoxNUU} \,\, Suppose that $N$ is a free group and that $T_s/T_{s+1}
$  and $N_{(s)}/ N_{(s+1)}$  are torsion-free for  $1\le s\le n$. Then there is a non-natural isomorphism
  \[ Q_n(G,N) \:\hcong\: \bigoplus_{i=0}^{n-1} \hspace{2mm}
{\rm U}_i{\rm L}^{\gamma} (T) \ot {\rm U}_{n-i-1}{\rm L}^{\cal N}(N) \ot N\ab
 \:.\]
\end{kor}

\section{Proofs for section 2}



%


The starting point of our approach is  the following elementary fact.

\begin{lem}\label{lemsplit}\quad Let $G$ be a group, $H$ and $K$ two subgroups of $G$ such that $H\cap K=\{1\}$, and let 
$J$\  be a left
ideal of $\Z(H)$ contained in $I(H)$. Then one has a short exact sequence
  \[   \IZ(K) \IZ(H) J \:\hra\: \IZ(K)  J \mr{s} \IZ(K)\ot  J /\IZ(H) J \,\to\, 0 \]
where $s((k-1)x) = (k-1) \ot(x + \IZ(H) J )$, $k\in K$, $x\in J$.
\end{lem}

\proof When the symbols $k$ resp.\ $h$ run through the nontrivial elements of $K$ resp.\ $H$ the elements $kh$ are distinct, and also distinct from the elements of $H$ and $K$. Thus the  map $\mu\mapco I(K) \ot I(H) \mr{} I(K)I(H)$ given by multiplication in $\Z(G)$ is an isomorphism since it sends the canonical basis $((k-1) \ot (h-1))$ of $I(K) \ot I(H)$ to linearly independant elements in $\Z(G)$.
Consider the following commutative square with $j\mapco J \hra I(H)$:
\[\begin{matrix}
\ruled I(K) J  & \hra  &  I(K) I(H) \cr
\surup{\mu'}  &  & \isoup{\mu} \cr
I(K) \ot J  & \mr{1\ot j}  &  I(K) \ot I(H)\ruled
\end{matrix}\]
The map $1\ot j$ is injective as $I(K)$ is a free \Z-module, hence $\mu'$ is an isomorphism, too. So we have the following commutative diagram with exact rows where $\mu''$ is given by restriction of $\mu'$, $\bar{\mu}'$ is induced by $\mu'$, $j'\mapco I(H)J \hra J$ is the injection and $q$ the corresponding quotient map.
\[\begin{matrix}
I(K)I(H)J  & \hra  &  I(K)J  &  \to  & I(K)J / I(K)I(H)J & \to & 0\ruled\cr
\surup{\mu''}  &  &  \isoup{\mu'}  &  &  \mapup{\bar{\mu}'} \cr
I(K)\ot I(H)J  & \mr{1\ot j'}   &  I(K) \ot J  &  \mr{1\ot q}  &  I(K) \ot J/ I(H)J  & \to & 0
\end{matrix}\]
This shows that $\bar{\mu}'$ is an isomorphism  which implies the assertion.\hfbox\V

In the sequel, we consider an {\em arbitrary}\/ descending filtration $\Delta\,\colon\: \IZ(N) = \Delta_1
\supset \Delta_2 \supset \cdots$ of $I(N)$ by subgroups $\Delta_i$; later on, we shall
specialize to the cases $\Delta  = \Lambda $ or $\Delta  = \IZ(N) \Lambda$ where $(\IZ(N) \Lambda)_i = \IZ(N)
\Lambda_{i-1}$. Let
  \[ \KnD{n} = \sum_{i=1}^{n-1} \Delta_{n-i}\,I^i(T) \:.\]
Then $\KnL{n} = {\cal K}_n$ while 
  \begin{eqnarray}
\KnINL{n}  &=&   \sum_{i=1}^{n-1} I(N)\,\Lambda_{n-i-1}\,I^i(T) \nonumber\\
  &=&  I(N)\,\sum_{i=1}^{n-1} \Lambda_{n-1-i}\,I^i(T) \nonumber\\
  &=&  \ruleu \ruled I(N)\,{\cal K}_{n-1}^* \:.\label{KnINL}
\end{eqnarray}

Thus computing filtration quotients $\KnD{n}/ \KnD{n+1}$ amounts to computing the direct factors ${\cal K}_n/{\cal K}_{n+1}$ and
$I(N){\cal K}_{n-1}^*/I(N){\cal K}_n^*$ of the abelian groups $\QnG{n} $, $\QnGT{n}$, and of $I(N)I^{n-1}(G)/I(N)I^{n}(G)$ instead of
$\QnGN{n}$, see \REF{AugTallg}, \REF{FoxTallg} and \REF{NG}. But the latter default is easily corrected by using the following device.\V

\begin{bem}\label{symbem} \rm 
  By Remark \ref{conjonFox} we obtain a commutatif diagram
\[\begin{matrix}
\frac{\dst I(N)I^{n-1}(G)}{\dst I(N)I^{n}(G)} & \mr{(-)^{\star}} & \frac{\dst I^{n-1}(G)I(N)}{\dst  I^{n}(G)I(N)\ruled}\cr
\| & & \|\cr
\frac{\dst \ruleu I(N)\Lambda_{n-1}}{\dst I(N) \Lambda_{n}} \hoplus \frac{\dst {\cal K}_n^{I(N)\Lambda}}{
\dst {\cal K}_{n+1}^{I(N)\Lambda}} 
& \isor{(-)^{\star}\oplus(-)^{\star}} & 
\frac{\dst \ruleu \Lambda_{n-1}I(N)}{\dst   \Lambda_{n}I(N)} \hoplus \frac{\dst \Gamma_{n-1}^* I(N)}{
\dst\Gamma_{n}^* I(N)} 

\end{matrix}\]

Our computation of the quotients $\KnD{n}/ \KnD{n+1}$ below provides a functorial computation of 
${\cal K}_{n}^{I(N)\Lambda}/{\cal K}_{n+1}^{I(N)\Lambda}$ in terms of induced and connecting maps between certain tensor and torsion products, namely between quotients of the filtrations $(I^i(N)\Lambda)_{i\ge 0}$ and $(I^j(T))_{j\ge 1}$ of $I(N)$ and $I(T)$, in this order. But it is easily checked that applying the symmetry isomorphisms of the tensor and torsion product, as well as the conjugation isomorphisms $(I(N)\Lambda)_i \hcong (\Lambda I(N))_i$, to our description is compatible with the conjugation isomorphism ${\cal K}_{n}^{I(N)\Lambda}/{\cal K}_{n+1}^{I(N)\Lambda} \hcong \Gamma_{n-1}^* I(N)/\Gamma_{n}^* I(N)$. Thus taking the ``mirror-symmetric" version of our description of the former quotients provides a description of the latter.
 So it finally suffices to determine the quotients $\Lambda_n/\Lambda_{n+1}$,
$I^n(T)/I^{n+1}(T)$, $\Lambda_{n-1}I(N)/\Lambda_{n}I(N)$, and
$\KnD{n}/\KnD{n+1}$ in order to determine  $I^{n-1}(G)I(K)/I^{n}(G)I(K)$ for $K=G,N$ and $T$.
\end{bem}


Let $1 \le i\le n-1$, $1\le j\le n-i$ and $i+1\le m\le \infty$. Putting $I^{\infty}(T) =0$ we have connecting homomorphisms $\tau_1=\tau_1^{n,i,j}$ and 
$\tau_2=\tau_2^m=\tau_1^{n,i,j,m}$
  \[ \frac{\dst \Delta_{n-i}}{\dst \Delta_{n-i+1}}  \ot \frac{\dst \IZ^i(T)}{\dst
\IZ^{i+1}(T)} 
\:\ml{\tau_1}\:
 {\rm Tor}_1^{\Z}\Big( 
\frac{\dst \Delta_j}{\dst \Delta_{n-i}}\,,\,
\frac{\dst \IZ^i(T)}{\dst \IZ^{i+1}(T)} \Big)  
\:\mr{\tau_2^m}\: 
\frac{\dst \Delta_j}{\dst \Delta_{n-i}} \ot \frac{\dst \IZ^{i+1}(T)}{I^m(T)} \]

\N obtained from the short exact sequences
   \BE \label{Deltasequ}    \Delta_{n-i}/  \Delta_{n-i+1}  \hspace{2mm}\:\hra\:\hspace{2mm} 
 \Delta_j/\Delta_{n-i+1}   \hspace{2mm}\to\hspace{2mm}  \Delta_j/\Delta_{n-i}\hspace{2mm}  \to 0 \EE
  \[  \IZ^{i+1}(T)/\IZ^{m}(T) \hspace{2mm}\hra\hspace{2mm} \IZ^{i }(T)/\IZ^{m}(T) \hspace{2mm}\to\hspace{2mm} \IZ^{i }(T)/\IZ^{i+1}(T)
\hspace{2mm}\to 0
\] Recall that for  a canonical generator $\langle
\bar{x} ,k,\bar{y} \rangle$  of $ {\rm Tor}_1^{\mathbb{Z}}(  \Delta_j / \Delta_{n-i}
\,,\,\IZ^{i }(T)/\IZ^{i+1}(T) )$, i.e. $x \in \Delta_j$, $y \in \IZ^{i }(T)$, $k \in \Z$ such
that $kx \in \Delta_{n-i}$ and $ky \in \IZ^{i+1}(T)$, one has
    \BE\label{tausplit} \begin{matrix} 
\tau_1 \langle \bar{x},k,\bar{y}\rangle & = & \overline{kx} \ot \bar{y} \cr
\tau_2 \langle \bar{x},k,\bar{y}\rangle & = & \bar{x} \ot \overline{ky}\,.\end{matrix} \EE\vspace{1mm}

\begin{prop}\label{amalgam}\quad For $1 \le i \le n-1$ there is an exact sequence 
 \[
 {\rm Tor}_1^{\mathbb{Z}}\Big( 
\frac{\dst \ruleu\Delta_1}{\dst \Delta_{n-i}}\,,\,
\frac{\dst \IZ^i(T)}{\dst \IZ^{i+1}(T)} \Big)  
\:\mr{(\mbox{} - \tau_1,\nu_i\tau_2^n)^t}  \: 
 \left(\frac{\dst \Delta_{n-i}}{\dst \Delta_{n-i+1}}  \hspace{1mm}\ot\hspace{1mm}  \frac{\dst \IZ^i(T)}{\dst
\IZ^{i+1}(T)} \right) 
\:\oplus\:
\frac{\dst \Delta_1 \IZ^{i+1}(T)}{\dst \sum_{k=i+1}^n \Delta_{n-k+1}\,\IZ^{k}(T)}  \]
 \BE\label{prototyp}  \hspace*{0cm} \hmr{(\mu_i,\iota_i)}   
\frac{\dst \Delta_1 \IZ^{i }(T)}{\dst \sum_{k=i}^n \Delta_{n-k+1}\,\IZ^{k}(T)} \hmr{s_i} \frac{\dst \Delta_1}{\dst \Delta_{n-i}} \hspace{1mm}\ot
\hspace{1mm} \frac{\dst \IZ^i(T)}{\dst  \IZ^{i+1}(T)} \hmr{} 0
 \EE
\ruleu where $\iota_i$ is induced by the injection $\Delta_1 \IZ^{i+1}(T) \hra \Delta_1 \IZ^{i }(T)$,  $s_i(\overline{xy}) = \bar{x}\ot \bar{y}$ for $(x,y) \in \Delta_1 \times \IZ^{i }(T)$, 
$\nu_i$ is given by
 \[\nu_i\,\colon\: \frac{\dst \Delta_1}{\dst \Delta_{n-i}}
 \hspace{1mm}\ot\hspace{1mm}  \frac{\dst \IZ^{i+1}(T)}{\dst
\IZ^{n}(T)} 
 \hspace{1mm}{\cong}\hspace{1mm}  
 \frac{\dst \Delta_1\ot  \IZ^{i+1}(T)}{\dst
{\rm Im}(\Delta_{n-i} \ot \IZ^{i+1}(T) + \Delta_1\ot 
\IZ^{n}(T))} \mr{\tilde{\nu}_i}
\frac{\dst \Delta_1 \IZ^{i+1}(T)}{\dst \sum_{k=i+1}^n \Delta_{n-k+1}\,\IZ^{k}(T)}  \]
with $\tilde{\nu}_i$ being induced by multiplication in $\Z(G)$, and $\mu_i$ is defined in a similar way.

\end{prop}\vspace{1.5mm}

\proof Consider the following diagram
  \[\begin{matrix} 
{\rm Tor}_1^{\mathbb{Z}}\Big( 
\frac{\ruleu \dst \Delta_1}{\dst \Delta_{n-i}\ruled}\,,\,
\frac{\dst \IZ^i(T)}{\dst \ruled\IZ^{i+1}(T)} \Big)  &
\Inj{\tilde{\tau}_1}  &  \Delta_{n-i} \ot \left(\frac{\dst \IZ^i(T)}{\dst \IZ^{i+1}(T)}\right)  &
\mr{\alpha \ot id}  &  \Delta_1 \ot \left(\frac{\dst \IZ^i(T)}{\dst \IZ^{i+1}(T)}\right)  \ruled\cr
\mapdown{\nu_i \tau_2^n}  &  &  \mapdown{\tilde{\mu}_i}  &  & \| \cr
\frac{\dst \ruleu\Delta_1 \IZ^{i+1}(T)}{\dst \sum_{k=i+1}^n \Delta_{n-k+1}\,\IZ^{k}(T)} &
\Inj{\tilde{\iota}_i}  & 
\frac{\dst \Delta_1 \IZ^{i }(T)}{\dst \sum_{k=i+1 \ruled}^n \Delta_{n-k+1}\,\IZ^{k}(T)}   &
\Sur{s}  &
\Delta_1 \ot \left(\frac{\dst \IZ^i(T)}{\dst \IZ^{i+1}(T) }\right)  \ruleu
\end{matrix}\]
The top row is part of a
six-term exact sequence associated with the short exact sequence $\Delta_{n-i}
\stackrel{\alpha}{\hra} \Delta_1 \auf \Delta_1/  \Delta_{n-i}$, and hence is exact; note that $\T{\tau}_1$ is injective since
$\Delta_1$ is a free \Z-module. 
The maps $\T{\mu}_i,\T{\iota}_i$ are given by multiplication and inclusion, repectively. The bottom row is induced by the exact sequence in Lemma \ref{lemsplit} for $(K,H,J)=(N,T,I^i(T))$ and hence is also exact. Moreover, the diagram commutes; to see this for the left-hand
square  use  \REF{tausplit}. Now an easy diagram chase together with right-exactness of the tensor product shows that the sequence 
  \[ 0 \hspace{1mm}\to\hspace{1mm} {\rm Tor}_1^{\mathbb{Z}}\Big( 
\frac{\dst \Delta_1}{\dst \Delta_{n-i}}
\,,\,
\frac{\dst \IZ^i(T)}{\dst \IZ^{i+1}(T)} \Big) 
\:\mr{(\mbox{} - \tilde{\tau}_1,\nu_i\tau_2^n)^t}  \: 
   \Delta_{n-i}  \ot \left(\frac{\dst \IZ^i(T)}{\dst
\IZ^{i+1}(T)} \right) 
\:\oplus\:
\frac{\dst \Delta_1 \IZ^{i+1}(T)}{\dst \sum_{k=i+1}^n \Delta_{n-k+1}\,\IZ^{k}(T)}\]
  \BE\label{preamalgam} \hspace*{1cm}
   \hmr{(\tilde{\mu}_i,\tilde{\iota}_i)}  
 \frac{\dst \Delta_1 \IZ^{i }(T)}{\dst \sum_{k=i+1}^n \Delta_{n-k+1}\,\IZ^{k}(T)} \hmr{(q_i\ot 1)s} \frac{\dst \Delta_1}{\dst \Delta_{n-i}} \hspace{1mm}\ot
\hspace{1mm} \frac{\dst \IZ^i(T)}{\dst  \IZ^{i+1}(T)} \hmr{} 0
   \EE
is  exact where $q_i\mapco \Delta_1 \auf \Delta_1/\Delta_{n-i}$ is the canonical projection. Then the assertion follows by passing to the quotient modulo ${\rm Im}\big(
\Delta_{n-i+1} \ot  (\IZ^i(T) / \IZ^{i+1}(T) )
 \big)$ and modulo $  \T{\mu}_i\,\mbox{Im}
\big(\Delta_{n-i+1} \ot  (\IZ^i(T) / \IZ^{i+1}(T)) 
 \big) = \Imm{\Delta_{n-i+1} \ot \IZ^i(T)} $, respectively. Just note that $\T{\tau}_1$ composed with the
quotient map  $\Delta_{n-i } \ot \left(\frac{ \IZ^i(T)}{ \IZ^{i+1}(T)}\right) \:\auf\: 
\left(\frac{  \Delta_{n-i}}{  \Delta_{n-i+1}}\right) \ot \left( \frac{  \IZ^i(T)}{ 
\IZ^{i+1}(T)}\right)$ equals $\tau_1$ by naturality of six-term exact sequences. \hfbox\vspace{4mm}

As we will see next, Proposition \ref{amalgam} allows to successively ``unscrew'' the filtration quotients of $\KnD{}$. The first
case, however, is plain:\vspace{4mm}

\N{\bf Computation of $\KnD{2}/\KnD{3}$:}\hspace{3mm} For $n=2$ and $i=1$ Proposition \ref{amalgam} provides the exact sequence
\[  0 \hspace{1mm} = \hspace{1mm} {\rm Tor}_1^{\mathbb Z}\Big( \frac{\Delta_1}{\Delta_1} \,, \frac{I(T)}{I^2(T)} \Big) \hmr{} 
\frac{\Delta_1}{\Delta_2} \hspace{2mm}\ot\hspace{2mm} \frac{I(T)}{I^2(T)}
\hspace{2mm}\oplus\hspace{2mm}   \frac{\Delta_1 \,I^2(T)}{\Delta_1 \,I^2(T)} \hspace{2mm}\mr{(\mu_1,\iota_1)}\hspace{2mm} 
 \frac{\KnD{2} }{ \KnD{3} } \hspace{1mm}\to\hspace{1mm} 0 \]
whence 
\begin{eqnarray}
\frac{\KnD{2} }{ \KnD{3} } &\cong&   \frac{\Delta_1}{\Delta_2} \hspace{2mm}\ot\hspace{2mm} \frac{I(T)}{I^2(T)} \nonumber\\
  &\cong& \label{K2/K3} \rule{0mm}{12mm} \left\{
\begin{array}{ccccl}
 \rule{0mm}{-7mm} N/N_{(2)} & \ot & T/T_2 & \mbox{\hspace{2mm} if \hspace{2mm}} & \Delta = \Lambda\\
N/N_{ 2 } & \ot & T/T_2 & \mbox{\hspace{2mm} if \hspace{2mm}} & \Delta = I(N)\Lambda
\end{array} 
\right.
\end{eqnarray}

as
\BE\label{INL1/INL2}  \frac{ (I(N)\,\Lambda)_1 }{ (I(N)\,\Lambda)_2 }  = \frac{I(N)}{I^2(N)} \hspace{2mm}\cong\hspace{2mm} N/N_2
\:.\EE
\V

\N \proofofthm{Fox2}  By \REF{AugTallg} and \REF{FoxTallg} the desired computation of $Q_2(G)$ and $Q_2(G,T)$ follows from
Theorem \ref{U1-U2} and \REF{K2/K3} for $\Delta=\Lambda$. By \REF{FoxNallg}, $\QnGN{2} \hspace{2mm}\cong\hspace{2mm} 
 \Lambda_1I(N)/\Lambda_2I(N)  \hspace{2mm}\oplus\hspace{2mm}   \Gamma_1^*I(N)/\Gamma_2^*I(N)$. But
$ \Lambda_1I(N)/\Lambda_2I(N)  \hspace{6pt}\cong\hspace{6pt} {\rm U}_2^{{\cal N}\gamma}(N,N)$ by Theorem \ref{thetaGH12}, and
 $\Gamma_1^*I(N)/\Gamma_2^*I(N) \hcong T/T_2 \ot N/N_2$ by Remark \ref{symbem} and \REF{K2/K3}. \hfbox\vspace{10pt}




\nc{\TOR}[2]{{\rm Tor}_1^{\mathbb Z}\Big(#1\,,#2\Big)}

\nc{\TOr}[2]{{\rm Tor}_1^{\mathbb Z}\big(#1\,,#2\big)}

\N{\bf Computation of $\KnD{3}/\KnD{4}$:}\hspace{3mm} Taking $n=3$ and $i=1,2$ Proposition \ref{amalgam} provides the following
two exact sequences
\[  \TOR{ \frac{\Delta_1}{\Delta_2} }{ \frac{I(T)}{I^2(T)} } \hmr{(-\tau_1,\nu_1\tau_2^3)^t} 
\frac{\Delta_2}{\Delta_3} \hspace{2mm}\ot\hspace{2mm} \frac{I(T)}{I^2(T)}
\hspace{2mm}\oplus\hspace{2mm}   \frac{\Delta_1 \,I^2(T)}{\Delta_2 \,I^2(T) + \Delta_1 \,I^3(T)}
\hspace{2mm}\Sur{(\mu_1,\iota_1)}\hspace{2mm} 
 \frac{\KnD{3} }{ \KnD{4} } \]
\[  \TOR{ \frac{\Delta_1}{\Delta_1} }{ \frac{I^2(T)}{I^3(T)} } \hmr{} 
\frac{\Delta_1}{\Delta_2} \hspace{2mm}\ot\hspace{2mm}  \frac{I^2(T)}{I^3(T)}
\hspace{2mm}\oplus\hspace{2mm}   \frac{\Delta_1 \,I^3(T)}{\Delta_1 \,I^3(T) }
\hspace{2mm}\Sur{(\mu_2,\iota_2)}\hspace{2mm} 
 \frac{\Delta_1 \,I^2(T)}{\Delta_2 \,I^2(T) + \Delta_1 \,I^3(T)} \ruled\]
\vspace{1mm}

Using Theorem \ref{U1-U2} we thus obtain the following result which requires the connecting homomorphisms 
\[ \frac{\Delta_{p+1}}{\Delta_{p+2}} \ot \UTn{q} \hspace{6pt}\ml{\hat{\tau}_1^{pq}} \hspace{6pt}
\TOr{\frac{\Delta_{p}}{\Delta_{p+1}}}{\UTn{q}} 
\hspace{6pt}\mr{\hat{\tau}_2^{pq}} \hspace{6pt} \frac{\Delta_{p }}{\Delta_{p+1}} \ot \frac{\UTn{q+1}}{\Imm{\delta_2^q}}\]
for $p,q=1,2$ where $\hat{\tau}_1^{pq}$ is induced by the short exact sequence \REF{Delp},
and $\hat{\tau}_2^{pq}$ is induced by sequence \REF{U2P2G}  for $q=1$ and by sequence \REF{U3U2sequ} for $q=2$, both with $(G,\calG)=(T,\gamma)$,
compare \REF{tau1tau2}.

\begin{satz}\label{KD3/KD4} For any descending subgroup filtration $\Delta$ of $I(N)$ there is a natural exact sequence
\[ \Tor{  \Delta_1/\Delta_2  }{ T\ab} \hmr{\delta_2^{\Delta}} 
 (\Delta_2/\Delta_3)  \hspace{0mm}\ot\hspace{0mm}   T\ab
 \hspace{2mm}\oplus\hspace{2mm}    (\Delta_1 /\Delta_2)   \hspace{0mm}\ot\hspace{0mm} \UTn{2}
\hspace{1mm}\Sur{\mu}\hspace{1mm} 
  \KnD{3} / \KnD{4}   \]

\N where  $\delta_2^{\Delta} = (-\hat{\tau}_1^{11}, \hat{\tau}_2^{11})^t$  and $\mu$ is given by the isomorphisms $\theta_k^{\gamma}$, $k=1,2$, and multiplication in $\Z(G)$.\hfbox
\end{satz}

\proofofthm{Fox3} By \REF{AugTallg} and \REF{FoxTallg} the desired computation of $\QnG{3}$ and $\QnGT{3}$ follows from Theorem 
 \ref{KD3/KD4} taking $\Delta=\Lambda$: it suffices to note that Theorem \ref{U1-U2} provides an isomorphism between  sequence \REF{Delp} for $p=1$ and sequence \REF{U2P2G} for $(G,\calG)=(N,\cal N) $; this isomorphism transforms the maps $\hat{\tau}_k^{11}$ in 
 $ {\tau}_k^{11}$, $k=1,2$, which were computed in \REF{tau1q} and \REF{taup1}.
 
Next by \REF{FoxNallg},  $\QnGN{3} \:\cong\: \Lambda_2 I(N)/\Lambda_3 I(N) \hspace{2mm}\oplus\hspace{2mm} 
\Gamma^*_2 I(N)/\Gamma^*_3 I(N)$. Taking $\Delta=I(N)\Lambda$ and using \REF{KnINL} Theorem \ref{KD3/KD4} provides a presentation
of $I(N){\cal K}_2^*/I(N){\cal K}_3^*$ which turns into the desired one of 
$\Gamma^*_2 I(N)/\Gamma^*_3 I(N)$ by means of  Remark \ref{symbem}. To make this explicit, first note that the mirror-symmetric version  of sequence \REF{Delp} for $p=1$ and $\Delta=I(N)\Lambda$ is
\BE\label{L2L3sequ}
0 \hspace{6pt}\to\hspace{6pt} \frac{\Lambda_1I(N)}{\Lambda_2I(N)} \hspace{6pt}\stackrel{\alpha}{\hra} \hspace{6pt} \frac{\Lambda_1}{\Lambda_2I(N)} \hspace{6pt}\to\hspace{6pt} \frac{\Lambda_1}{\Lambda_1I(N)} \hspace{6pt}\to\hspace{6pt} 0
\EE  
Next, the mirror-symmetric versions  of the maps $\hat{\tau}_1^{11}$ and $\hat{\tau}_2^{11}$ for $\Delta=I(N)\Lambda$ are the connecting homomorphisms $\tau_2^{\Lambda I(N)}$ and 
$\tau_1^{\Lambda I(N)}$, resp., which form the top row of the commutative diagram
\[\begin{matrix}
T\ab \ot \frac{\dst \Lambda_1I(N)}{\dst \Lambda_2I(N)}  & \ml{\tau_2^{\Lambda I(N)}} & \TOr{T\ab}{\frac{\dst \Lambda_1}{\dst \ruled\Lambda_1I(N)}} & \mr{\tau_1^{\Lambda I(N)}} & \UTn{2} \ot \frac{\dst \Lambda_1}{\dst \Lambda_1I(N)} \cr
\isoup{1\ot \theta_2^{\cal N\gamma}(\overline{\nu_{11}^{\cal N\gamma}})} & & \isoup{(1,\theta_1^{\cal N\gamma})_*} & & \isoup{1\ot \theta_1^{\cal N\gamma}}\cr
T\ab \ot \frac{\dst N\AB\ot N\ab}{\dst (q_2^{\cal N}\ot 1)l_2^{\gamma}\,{\rm Ker}(c_2^{\gamma})} &
\ml{\tilde{\tau}_2^{\Lambda I(N)}} &
\Tor{T\ab}{N\ab} &
\mr{\tilde{\tau}_1^{\Lambda I(N)}} &
\UTn{2} \ot N\ab

\end{matrix}
\]
and which are induced by the sequences \REF{L2L3sequ}, and \REF{U2P2G} for $(G,\calG)=(T,\gamma)$, resp. 

To compute $\tilde{\tau}_2^{\Lambda I(N)}$ let $\langle tT_2,k,nN_2\rangle$ be a canonical generator of $\Tor{T\ab}{N\ab}$. Then
\begin{eqnarray*}
\tilde{\tau}_2^{\Lambda I(N)}\langle tT_2,k,nN_2\rangle &=&  (1\ot \theta_2^{\cal N\gamma}(\overline{\nu_{11}^{\cal N\gamma}}))^{-1} \tau_2^{\Lambda I(N)} (1,\theta_1^{\cal N\gamma})_*\langle tT_2,k,nN_2\rangle \\
&=& (1\ot \theta_2^{\cal N\gamma}(\overline{\nu_{11}^{\cal N\gamma}}))^{-1} \tau_2^{\Lambda I(N)}\langle tT_2\,,k\,,\,n-1+I^2(N) \rangle \\
&=& (1\ot \theta_2^{\cal N\gamma}(\overline{\nu_{11}^{\cal N\gamma}}))^{-1} \Big( (tT_2) \ot \alpha^{-1}\big( k(n-1)+\Lambda_2I(N)\big)\Big)\\
&=& (tT_2) \ot  (\theta_2^{\cal N\gamma}(\overline{\nu_{11}^{\cal N\gamma}}))^{-1}   \Big( (n^k-1)-{k\choose 2}(n-1)^2+\Lambda_2I(N)\Big)\\
&=& (tT_2) \ot  (\theta_2^{\cal N\gamma}(\overline{\nu_{11}^{\cal N\gamma}}))^{-1}   \Big( (\prod_{i=1}^p [n_i,n_i\st] -1)-{k\choose 2}(n-1)^2+\Lambda_2I(N)\Big)\\
&=& (tT_2) \ot  (\theta_2^{\cal N\gamma}(\overline{\nu_{11}^{\cal N\gamma}}))^{-1}   \Big( \sum_{i=1}^p  \Big((n_i-1)(n_i\st-1) - (n\st_i-1)(n_i-1)\Big)\\
& & \hspace{55pt} {}-{k\choose 2}(n-1)^2+\Lambda_2I(N) \Big)\quad\mbox{by \REF{ab-1} and \REF{comm}}\\
&=& (tT_2) \ot   \bigg( \sum_{i=1}^p  \Big((n_iN_{(2)})\ot (n_i\st N_2) - (n\st_iN_{(2)})\ot (n_iN_2)\Big)\\
& & \hspace{50pt}{}-{k\choose 2}(nN_{(2)})\ot (nN_{2}) + (q_2^{\cal N}\ot 1)l_2^{\gamma}\,\Ker{c_2^{\gamma}}
\bigg)\\
\end{eqnarray*}


This provides  the first component of the map $\delta_3=(\tilde{\tau}_2^{\Lambda I(N)},{}-\tilde{\tau}_1^{\Lambda I(N)})$ which is the mirror-symmetric version of ${}-\delta_2^{I(N)\Lambda }$; the second component is obtained by a similar computation of $\tilde{\tau}_1^{\Lambda I(N)}$.\hfbox\vspace{15pt}

\nc{\SC}{\scriptstyle}

\N{\bf Computation of $\KnD{4}/\KnD{5}$:}\hspace{3mm} Proposition \ref{amalgam} provides the following four
 exact sequences extracted from sequence \REF{prototyp}, taking $n=4$ and $i=1,2,2,3$ resp.: 
\BE\label{n=4,i=1}  \begin{matrix} \Tor{ \frac{\Delta_1}{ \Delta_3} }{ \frac{I(T)}{I^2(T)} } \hmr{(-\tau_1,\nu_1\tau_2^4)} 
\frac{\Delta_3}{\Delta_4} \hspace{0.8mm}\ot\hspace{0.8mm} \frac{I(T)}{I^2(T)}
\hspace{2mm}\oplus\hspace{2mm}   \frac{\Delta_1 \,I^2(T)}{\sum_{k=2}^4 \Delta_{5-k}I^k(T)} \hspace{1mm}\Sur{(\mu_1,\iota_1)} \hspace{1mm} \frac{\Delta_3 \,I(T) +\Delta_1 \,I^2(T)}{\KnD{5}\ruled}\end{matrix}
 \EE
\BE\label{n=4,i=2,a} \begin{matrix} \Tor{ \ruleu\frac{\Delta_1}{\Delta_2} }{ \frac{I^2(T)}{I^3(T)} }  \hmr{(-\tau_1,\nu_2\tau_2^4)}  
\frac{\Delta_2}{\Delta_3} \hspace{0.8mm}\ot\hspace{0.8mm}  \frac{I^2(T)}{I^3(T)}
\hspace{2mm}\oplus\hspace{2mm}   \frac{\Delta_1 \,I^3(T)}{\sum_{k=3}^4 \Delta_{5-k}I^k(T) }
\hspace{2mm}\Sur{(\mu_2,\iota_2)}\hspace{2mm} 
 \frac{\Delta_2 \,I^2(T) + \Delta_1 \,I^3(T)}{\sum_{k=2}^4 \Delta_{5-k}I^k(T)\ruled}  \end{matrix}\EE
\BE\label{n=4,i=2,b} \begin{matrix}   \frac{\ruleu \Delta_2 \,I^2(T) + \Delta_1 \,I^3(T)}{\sum_{k=2}^4 \Delta_{5-k}I^k(T)} \hspace{2mm} \hra \hspace{2mm} 
 \frac{\Delta_1 \,I^2(T)  }{\sum_{k=2}^4 \Delta_{5-k}I^k(T)} \hspace{2mm} \Sur{s_2}\hspace{2mm}  \frac{\Delta_1}{ \Delta_2} \hspace{0.8mm}\ot\hspace{0.8mm}  \frac{I^2(T)}{I^3(T)} 
 \ruled \end{matrix}\EE
 \BE\label{n=4,i=3}  \begin{matrix} \Tor{ \frac{\ruleu \Delta_1}{ \Delta_1} }{ \frac{I^3(T)}{I^4(T)} } \hmr{} 
\frac{\Delta_1}{\Delta_2} \hspace{0.8mm}\ot\hspace{0.8mm} \frac{I^3(T)}{I^4(T)}
\hspace{2mm}\oplus\hspace{2mm}   \frac{\Delta_1 \,I^4(T)}{\Delta_1 \,I^4(T)} \hspace{2mm}\Sur{(\mu_3,\iota_3)}\hspace{2mm}   \frac{\Delta_1 \,I^3(T)}{\sum_{k=3}^4 \Delta_{5-k}I^k(T) } \end{matrix}
 \EE
Combining  sequences \REF{n=4,i=2,a} and \REF{n=4,i=3} with Theorems \ref{U1-U2} and \ref{Q3} provides an exact sequence
 \BE\label{n=4,i=2,3}  \begin{matrix} 
 \Tor{ \frac{\ruleu \Delta_1}{ \Delta_2} }{ \UTn{2} } \hmr{(-\hat{\tau}_1^{12},\hat{\tau}_2^{12})}

 \frac{\ruleu \Delta_2}{ \Delta_3} \hspace{0.3mm}\ot\hspace{0.3mm}  
 \UTn{2} \hspace{1mm}\oplus\hspace{1mm}   \frac{\ruleu \Delta_1}{ \Delta_2}  \hspace{0.3mm}\ot\hspace{0.3mm} \coker{\delta_1^{\gamma}}  \hspace{0.8mm}\Sur{(\mu_2,\iota_2)} \hspace{0.8mm} \frac{\Delta_2 \,I^2(T) + \Delta_1 \,I^3(T)}{\sum_{k=2}^4 \Delta_{5-k}I^k(T)\ruled}  \end{matrix}\EE

\vspace{1mm}

\begin{satz}\label{KD4/KD5} For any descending subgroup filtration $\Delta$ of $I(N)$ the quotient $\KnD{4}/\KnD{5}$ is determined by the following tower of successive natural quotients where $\Ker{\pi_k} = \Imm{\xi_k}$, $k=1,2,3$.

\N\makebox[14.7cm]{ \makebox[0mm]{
\begin{minipage}{17cm}\small

\[ \begin{matrix}
\ruled\frac{\Delta_1}{\Delta_2} \ot \Tor{T\ab}{T\ab}  &  \mr{\xi_1}  &  \frac{\Delta_3}{\Delta_4} \ot T\ab \hoplus \frac{\Delta_2}{\Delta_3} \ot \UTn{2} 
 \hoplus \frac{\Delta_1}{\Delta_2} \ot \UTn{3} \cr
  & & \surdown{\pi_1} \cr
\ruleu\ruled\Tor{\frac{\Delta_2}{\Delta_3}}{T\ab} \oplus \Tor{\frac{\Delta_1}{\Delta_2}}{\UTn{2}}  & \mr{\xi_2}  &  \frac{\Delta_3}{\Delta_4} \ot T\ab \hoplus \frac{\Delta_2}{\Delta_3} \ot \UTn{2} 
 \hoplus \frac{\Delta_1}{\Delta_2} \ot \coker{\delta_1^{\gamma}} \cr
   & & \surdown{\pi_2} \cr
\ruleu\ruled \Ker{\delta_2^{\Delta}} & \mr{\xi_3}  &  \coker{\xi_2} \cr
      & & \surdown{\pi_3} \cr
  & & \ruleu \KnD{4}/\KnD{5}
\end{matrix}\]

\end{minipage}\ruled
}\rule[-11mm]{0mm}{3mm} }\vspace{6mm}
Here $\xi_1=(0,0,1\ot \delta_1^{\gamma})^t$, $\xi_2 = \left(\begin{matrix} {} -\hat{\tau}_1^{21} & \hat{\tau}_2^{21} & 0 \cr 0 & {} -\hat{\tau}_1^{12} & \hat{\tau}_2^{12}\end{matrix}\right)^t$, 
$\delta_2^{\Delta} = (-\hat{\tau}_1^{11}, \hat{\tau}_2^{11})^t$ as in Theorem \ref{KD3/KD4},
and the map $\xi_3$ is defined in the proof below.
To describe it explicitly, suppose that $
\Delta_1/\Delta_2$ and $T$ are finitely generated with cyclic decomposition ${\rm Tor}(\Delta_1/\Delta_2) = \bigoplus_{i=1}^r \Z/a_i\Z \langle \bar{x}_i \rangle$, $x_i\in \Delta_1$, and with the remaining notations of Theorem \ref{Fox4}. Then an element $\omega = \sum_{i,j} \langle \bar{x}_i, k_{ij},\bar{t}_j \rangle \in\Tor{\frac{\Delta_1}{\Delta_2}}{T\ab}$ lies in $\Ker{\delta_2^{\Delta}}$ if and only if the conditions (i) and (ii) in Theorem \ref{Fox4} hold, as well as the following condition\V

  (iii)'\hspace{2mm} $\forall 1\le j\le s$, \hspace{5mm} $\sum_{i=1}^r k_{ij} x_i = b_j\delta_j^2 + \delta_j^3$ with $\delta_j^k \in \Delta_k$, $k=2,3$.\V

   \N In this case, \begin{eqnarray*}
  \xi_3(\omega)  &=& \pi_2\pi_1 \left( {} - \sum_{j=1}^s \overline{\delta_j^3} \ot \bar{t}_j \hspace{1mm}\mathbb{,}\hspace{1mm} \sum_{i=1}^r  \overline{a_ix_i} \ot \bigg( (u_iT_3) -  \sum_{j=1}^s \frac{p_{ij}}{d_{ij}} {k_{ij} \choose 2} (t_jT_2)^2 \bigg) \right.\\
   & & {} - \sum_{j=1}^s \overline{\delta_j^2} \ot \bigg( (t_j^{b_j}T_3) - {b_j \choose 2} (t_jT_2)^2 \bigg) \hspace{1mm}\mathbb{,}\hspace{1mm} \\
  & &  \left.{} + \sum_{i=1}^r \bar{x}_i \ot  
  \left(\rule{0mm}{20pt}
  ( \bar{x}_i \ot   (v_iT_4) - \sum_{j=1}^s \left( \frac{q_{ij}}{d_{ij}} {k_{ij} \choose 2}    (t_j^{b_j}T_3)(t_jT_2)\right.\right.\right.\\
&& \left. \left.\left. {} +\left( {k_{ij} \choose 3}
   -   \frac{q_{ij}}{d_{ij}} {k_{ij} \choose 2} {b_j \choose 2}
    \right) (t_jT_2)^3 
\right) \rule{0mm}{20pt}\right)\rule{0mm}{23pt}\right)
    \end{eqnarray*}

\end{satz}

\proof First note that \vspace{-1mm}

\Ph\hfill$\coker{\xi_1}\hspace{1mm}\cong\hspace{1mm} \frac{\Delta_3}{\Delta_4} \ot T\ab \hoplus \frac{\Delta_2}{\Delta_3} \ot \UTn{2} 
 \hoplus \frac{\Delta_1}{\Delta_2} \ot \frac{I^3(T)}{I^4(T)}$\hfill\V
 
 \N by Theorem \ref{Q3}. Now consider the following diagram where implicitly $n=4$, $\alpha$ and $\beta$ are the obvious injections, $\alpha'=\alpha\ot 1$, $\rho,\rho',q$ are the obvious quotient maps, 
 and $\tau_1=\tau_1^{4,1,1}$, 
$\tau_2=\tau_1^{4,1,1,4}$, 
$\tau_1'=\tau_1^{4,1,2}$, 
$\tau_2'=\tau_1^{4,1,2,4}$.
 
 
 \N\makebox[14.7cm]{ \makebox[0mm]{
\begin{minipage}{17cm}\small
\BE\label{KD4/KD5dia}
\begin{matrix}
\Tor{\frac{\ruleu\Delta_2}{\ruled\Delta_3}}{ \frac{I(T)}{I^2(T)}}  &  \mr{(\alpha,1)_*}  &  \Tor{\frac{\Delta_1}{\Delta_3}}{ \frac{I(T)}{I^2(T)}}  &  \mr{(\rho,1)_*}  &  \Tor{\frac{\Delta_1}{\Delta_2}}{ \frac{I(T)}{I^2(T)}} \cr
\mapdown{(-\tau_1',\tau_2')^t}  &  &  \mapdown{(-\tau_1,\tau_2)^t}  &  &  \mapdown{\tau_2^{3,1,1,3}}  \cr
\frac{\ruleu\Delta_3}{\Delta_4} \ot \frac{I(T)}{I^2(T)}  \hhoplus \frac{\Delta_2}{\Delta_3} \ot \frac{I^2(T)}{I^4(T)}  &  \stackrel{1 \oplus \alpha'}{\hra}  & \frac{\ruleu\Delta_3}{\Delta_4} \ot \frac{I(T)}{I^2(T)}  \hhoplus
\frac{\Delta_1}{\Delta_3} \ot \frac{I^2(T)}{I^4(T)}  &
\Sur{(0,{\rho}\ot \rho')} & \frac{\Delta_1}{\Delta_2} \ot \frac{I^2(T)}{I^3(T)\ruled} \cr
\mapdown{1\oplus\nu_1(i\ot1)} & & \mapdown{1\oplus\nu_1} & & \| \cr
\frac{\ruleu\Delta_3}{\Delta_4} \ot \frac{I(T)}{I^2(T)}  \hhoplus \frac{\Delta_2 I^2(T) + \Delta_1 I^3(T)}{\sum_{k=2}^4 \Delta_{5-k}I^k(T)}  &  \stackrel{1 \oplus \beta}{\hra}  &
\frac{\Delta_3}{\Delta_4} \ot \frac{I(T)}{I^2(T)}  \hhoplus \frac{ \Delta_1 I^2(T)}{\sum_{k=2}^4 \Delta_{5-k}I^k(T)} & \Sur{(0,s_2)} & \frac{\Delta_1}{\Delta_2} \ot \frac{I^2(T)}{I^3(T)\ruled} \cr
\surdown{q}  &  & \surdown{(\mu_1,\iota_1)}  &  & \cr
\coker{(-\tau_1',\nu_1(i\ot 1)\tau_2')^t}  &  \mr{\overline{1 \oplus \beta}}  & \frac{\ruleu\Delta_3I(T) + \Delta_1 I^2(T)}{\KnD{5}}  &  & 
\end{matrix}
\EE


\end{minipage}\ruled
}\rule[-11mm]{0mm}{3mm} }\vspace{6mm}

The upper squares commute by naturality of the connecting homomorphisms; for the right-hand square we have to use this argument twice, according to the factorization $\rho \ot \rho' = (1\ot \rho')(\rho \ot 1)$: we have $(\rho \ot 1)\tau_2=\tau_2^{3,1,1,4} (\rho,1)_*$ and $(1\ot \rho')\tau_2^{3,1,1,4} = \tau_2^{3,1,1,3} $. The middle squares commute by definition of the concerned maps; together with exactness of sequence \REF{n=4,i=1} this implies that the map
$\overline{1\oplus\beta}$  induced by $1\oplus\beta$  is well-defined. 
So the whole diagram commutes. Furthermore,
the first and third row are exact; in fact, the top row is part of a 6-term exact sequence, and exactness of the third row  follows from sequence  \REF{n=4,i=2,b}. 

Now observe that $\KnD{4}/\KnD{5} = \Imm{\overline{1 \oplus \beta}}$. But using sequence \REF{n=4,i=2,3} and Theorems \ref{U1-U2} and \ref{Q3} we obtain an isomorphism
\[ \Xi \mapco \coker{(-\tau_1',\nu_1(i\ot 1)\tau_2')^t} \hspace{1mm}\isor{}\hspace{1mm} \coker{\xi_2}\,.\]
 Thus in order to establish the desired tower it remains to show that $\Xi ( \Ker{\overline{1 \oplus \beta}}) $ $= \Imm{\xi_3}$ for a suitable map $\xi_3$. To do this we apply the snake lemma with respect to the diagram whose rows are the first and third row of diagram \REF{KD4/KD5dia} and whose vertical maps are the compositions of the  respective vertical maps in \REF{KD4/KD5dia}.  Using sequence \REF{n=4,i=1} we thus obtain an exact sequence \newline
 \Ph\hfill$  \Imm{(\rho,1)_*} \cap \Ker{\tau_2^{3,1,1,3}} \mr{\xi_3'} \coker{(-\tau_1',\nu_1(i\ot 1)\tau_2')^t} 
\hspace{1mm} \mr{\overline{1 \oplus \beta} }\hspace{1mm}
 \frac{\Delta_3I(T) + \Delta_1 I^2(T)}{\KnD{5}}$\hfill
 \medskip

\N where $\xi_3'$  is given by the switchback rule $\xi_3' =  q (1\oplus \beta)^{-1} (-\tau_1, \nu_1\tau_2)^t (\rho,1)_*^{-1}$.
But $\Ker{\tau_2^{3,1,1,3}} = 
(1,\theta_1^{\gamma})_*\Ker{\hat{\tau}_2^{11}}$ 
and $\Imm{(\rho,1)_*}=\Ker{\tau_1^{3,1,1}} = (1,\theta_1^{\gamma})_*\Ker{\hat{\tau}_1^{11}}$ by the corresponding 6-term exact sequence, so we get $\Imm{(\rho,1)_*} \cap \Ker{\tau_2^{3,1,1,3}} =$ $ (1,\theta_1^{\gamma})_*\Ker{\delta_2^{\Delta}}$ since $(1,\theta_1^{\gamma})_*$ is an isomorphism. Now define $\xi_3$ to be 
$\Xi \xi_3' $ precomposed with the isomorphism $\Ker{\delta_2^{\Delta}} \isor{} \Imm{(\rho,1)_*} \cap \Ker{\tau_2^{3,1,1,3}}$ induced by $(1,\theta_1^{\gamma})_*$. Then the asserted tower is established.


 It remains to prove the claims concerning the element $\omega$. 
One has
 \[ \hat{\tau}_1^{11}(\omega) = \sum_{j=1}^s \Big( \sum_{i=1}^r k_{ij}x_i + \Delta_3 \Big) \ot (t_jT_2) \hspace{2mm}\in \frac{\Delta_2}{\Delta_3} \ot T\ab
 \hspace{2mm}\cong\hspace{2mm}
   \bigoplus_{j=1}^s \frac{\Delta_2}{\Delta_3} \ot \Z/b_j\Z \langle t_jT_2 \rangle\,,\]
  whence $\hat{\tau}_1^{11}(\omega) =0$  iff $\forall j$, $\sum_{i=1}^r k_{ij}x_i + \Delta_3 \in b_j \frac{\Delta_2}{\Delta_3}$, i.e., 
$\sum_{i=1}^r k_{ij}x_i = b_j \delta_j^2 + \delta_j^3$ with $\delta_j^k \in \Delta_k$. On the other hand, noting $ \bar{x_i}  =x_i+\Delta_2$ and $b_{jk} ={\rm gcd}(b_j,b_k)$, we obtain as in the proof of Lemma \ref{del1comp}
\begin{eqnarray*}
\hat{\tau}_2^{11}(\omega) &=&  \sum_{i=1}^r \bar{x_i} \ot 
\sum_{j=1}^s (\theta_2^{\gamma})^{ -1}\Big((k_{ij}(t_j-1)+ I^3(T)\Big) \\
  &=&  
\sum_{i=1}^r \bar{x_i} \ot  \sum_{j=1}^s \Big( (t_j^{k_{ij}}T_3) -  {k_{ij} \choose 2} (t_jT_2)^2 \Big) \\
&=& \sum_{i=1}^r \bar{x_i} \ot  \bigg(\prod_{j=1}^s   t_j^{k_{ij}} \bigg)T_3   \hspace{2mm} -  \hspace{2mm}
\sum_{i,j} \bar{x_i} \ot    {k_{ij} \choose 2} (t_jT_2)^2   \\
&\in& \frac{\Delta_1}{\Delta_2} \ot \UTn{2}\\
&\cong& 
\bigoplus_{i=1}^r \Big( \Z/a_i\Z \langle \bar{x}_i \rangle \ot (T_2/T_3) \Big) \oplus \bigoplus_{
\begin{matrix}
\scriptstyle 1\le i\le r\cr
\scriptstyle  1\le j\le k\le s
\end{matrix}} \Z/a_i\Z \langle \bar{x}_i \rangle \ot \Z/b_{jk}\Z \langle ( t_jT_2)\ot ( t_kT_2) \rangle 
\end{eqnarray*}
where we use the classical direct sum decomposition $\UTn{2} = T_2/T_3 \oplus {\rm SP}^2(T\ab)$ for finitely generated $T\ab$, see also Proposition \ref{U2NNsequ}. Thus $\hat{\tau}_2^{11}(\omega) =0$ iff $\forall i$, $\prod_{j=1}^s   t_j^{k_{ij}} \in T_2^{a_i}T_3$ and $\forall i,j$, $d_{ij}=\gcd(a_i,b_{jj})$ divides ${k_{ij} \choose 2}$.

But $d_{ij}$ divides $k_{ij}$ since the latter is a commun multiple of $a_i$ and $b_j$, as we suppose that the element $\langle \overline{x_i},k_{ij},\overline{t_j})\rangle$ is defined. Thus $d_{ij}$ divides ${k_{ij} \choose 2}$ if and only if condition (i) is satisfied.
It follows that $\omega\in \Ker{\delta_2^{\Delta}}$ iff the conditions (i), (ii) and (iii)' are satisfied. 

To simplify the notation we henceforth identify $T\ab$ with $I(T)/I^2(T)$ via the isomorphism $\theta_1^{\gamma}$. Then for $\omega\in \Ker{\delta_2^{\Delta}}$,
\[ \omega  =  \sum_{j=1}^s  \sum_{i=1}^r \langle \frac{k_{ij}}{b_j} \bar{x}_i  \,, b_j,t_jT_2 \rangle  =  
 \sum_{j=1}^s   \langle \sum_{i=1}^r \frac{k_{ij}}{b_j} \bar{x}_i\,, b_j,t_jT_2 \rangle  =  (\rho,1)_*(\omega')\] with 
$\omega' = \sum_{j=1}^s  \langle \sum_{i=1}^r \frac{k_{ij}}{b_j} x_i - \delta_j^2 + \Delta_3\,, b_j,t_jT_2 \rangle$. Abbreviating $A=\sum_{k=2}^4 \Delta_{5-k}I^k(T)$, we next get
\begin{eqnarray}
 & &\hspace{-10mm}(-\tau_1,\nu_1\tau_2)^t(\omega')  \nonumber \\
 &=&  \sum_{j=1}^s \bigg( {}- (\delta_j^3 + \Delta_4) \ot (t_jT_2) \hspace{1mm},\hspace{1mm}
\Big( \sum_{i=1}^r \frac{k_{ij}}{b_j} x_i - \delta_j^2 \Big) b_j(t_j-1) + A \bigg) \nonumber \\
&=&\label{tauomegast}  \left(  \hspace{-1mm} - \sum_{j=1}^s  \overline{\delta_j^3}   \ot \bar{t}_j 
\hspace{1mm},\hspace{1mm}
 \sum_{i=1}^r x_i \bigg( \sum_{j=1}^s  k_{ij}(t_j-1) \bigg)  -   \sum_{j=1}^s    \delta_j^2 \left( (t_j^{b_j} - 1) - {b_j \choose 2} (t_j-1)^2 \right) +A  \right) \nonumber \\
\end{eqnarray}

The remaining calculations are again based on the identities \REF{binom} and \REF{ab-1}.
Using these and  the fact that $t_j^{k_{ij}} \in T_2$ we obtain the following congruences modulo $I^4(T)$, for $1\le i\le r$:
\begin{eqnarray*}
 & &\hspace{-10mm}\sum_{j=1}^s  k_{ij}(t_j-1)  \\
 &\equiv&   \sum_{j=1}^s  \Big(  (t_j^{k_{ij}} - 1) - {k_{ij} \choose 2} (t_j-1)^2 - {k_{ij} \choose 3} (t_j-1)^3 \Big)\\
  &\makebox[0mm]{$\equiv$}\makebox[0mm]{\raisebox{-13pt}{\small (i)}}
  &  \prod_{j=1}^s t_j^{k_{ij}} - 1 - \sum_{j=1}^s \left( \frac{1}{d_{ij}} {k_{ij} \choose 2} \Big( a_i p_{ij} (t_j-1)^2 + b_j q_{ij} (t_j-1)^2 \Big)  +  {k_{ij} \choose 3} (t_j-1)^3 \right) \\
    &\makebox[0mm]{$\equiv$}\makebox[0mm]{\raisebox{-13pt}{\small (ii)}}
  &
    u_i^{a_i}v_i  - 1 - a_i \sum_{j=1}^s \frac{p_{ij}}{d_{ij}} {k_{ij} \choose 2}  (t_j-1)^2 
        - \sum_{j=1}^s  \left( \frac{q_{ij}}{d_{ij}} {k_{ij} \choose 2}  \Big( (t_j^{b_j} - 1) -  {b_j \choose 2}(t_j-1)^2  \Big) \right.\\
        & & \left. \times (t_j-1) +  {k_{ij} \choose 3} (t_j-1)^3  \right) \\
        \end{eqnarray*}
As $u_i^{a_i}v_i  - 1 \equiv {a_i} (u_i - 1) + (v_i  - 1)$ modulo $I^4(T)$ since $u_i, v_i \in T_2$ we obtain 
\begin{eqnarray} x_i \sum_{j=1}^s  k_{ij}(t_j-1) 
    &\equiv& a_ix_i \Big( (u_i - 1) - \sum_{j=1}^s \frac{p_{ij}}{d_{ij}} {k_{ij} \choose 2}  (t_j-1)^2 \Big)   \nonumber \\ 
& &{}     + 
    x_i \Bigg( (v_i  - 1) 
    - \sum_{j=1}^s  \left( \frac{q_{ij}}{d_{ij}} {k_{ij} \choose 2}   (t_j^{b_j} - 1)(t_j-1)  \right. \nonumber \\ & & {}+\left( {k_{ij} \choose 3}  -  \frac{q_{ij}}{d_{ij}} {k_{ij} \choose 2} {b_j \choose 2}\right) (t_j-1)^3  \bigg)  \Bigg) \hspace{6pt}\mbox{mod $A$} \label{xisum}
 \end{eqnarray}
The formula for $\xi_3(\omega)$ now follows by combining the identities \REF{tauomegast} and \REF{xisum}.\hfbox\vspace{5mm}

\proofof{Theorem \ref{Fox4}} Taking $\Delta= \Lambda$ the tower in Theorem \ref{KD4/KD5} transforms into the one in Theorem \ref{Fox4} by use of Theorems \ref{U1-U2} and \ref{Q3}. It remains to check equivalence between the set of conditions (i) - (iii) and the set (i), (ii) and (iii)', as well as the asserted formula for $\xi_3$. First of all, we may take $x_i=n_i-1$ in Theorem \ref{KD4/KD5}. Now let $\omega = \sum_{i,j} \langle n_i N_{(2)}, k_{ij}, t_jT_2 \rangle \in\Tor{N\AB}{T\ab} $. Using \REF{binom}  and \REF{ab-1} together with the fact that $n_i^{k_{ij}} \in N_{(2)}$ we get the following identities for $1\le j\le s$, denoting $a_{ij}={\rm gcd}(a_i,a_j)$:
\begin{eqnarray*}
(\theta_2^{\cal  N})^{-1} \Big( \sum_{i=1}^r k_{ij} (n_i-1) + \Lambda_3 \Big) &=& (\theta_2^{\cal  N})^{-1} \bigg( \sum_{i=1}^r \Big( (n_i^{k_{ij}} - 1) - {k_{ij} \choose 2} (n_i-1)^2 + \Lambda_3\Big) \bigg) \\
  &=&   \bigg(\prod_{i=1}^r  n_i^{k_{ij}} \bigg)N_{(3)}   - \sum_{i=1}^r {k_{ij} \choose 2} (n_iN_{(2)})^2   \\
 &\in& {\rm U}_2{\rm L}^{\cal N}(N) \\
    &=&  N_{(2)}/N_{(3)}  \oplus \bigoplus_{1\le i\le j\le r}  \Z/a_{ij}\Z \langle (n_iN_{(2)})(n_jN_{(2)}) \rangle
\end{eqnarray*}
Here we use the standard direct sum decomposition $\UNn{2} \hcong N_{(2)}/N_{(3)}  \hoplus {\rm SP}^2(N\AB)$, cf.\ \REF{U2SP2}. Hence (iii)' implies (iii). Conversely, suppose that $\omega$ satisfies conditions (i) - (iii). Then a similar calculation as in the proof of Theorem \ref{KD4/KD5} shows that for $1\le j\le s$ and modulo $\Lambda_4$,
\begin{eqnarray*}  \sum_{i=1}^r  k_{ij}(n_i - 1) 
    &\equiv& b_j \bigg( (y_j - 1) - \sum_{i=1}^r \frac{q_{ij}}{d_{ij}} {k_{ij} \choose 2}  (n_i-1)^2 \bigg)  
     + 
     (z_j  - 1)  \nonumber \\ 
    & & - \sum_{i=1}^r  \left( \frac{p_{ij}}{d_{ij}} {k_{ij} \choose 2}   (n_i^{a_i} - 1)(n_i-1) 
     \right. \nonumber \\ 
     & & {}+\left( {k_{ij} \choose 3}  -  \frac{p_{ij}}{d_{ij}} {k_{ij} \choose 2} {a_i \choose 2}\right) (n_i-1)^3  \bigg)     \nonumber \\
 \end{eqnarray*}

\N Thus condition (iii)' is satisfied for 
 \[
  \delta_j^2 =   (y_j - 1) - \sum_{i=1}^r \frac{q_{ij}}{d_{ij}} {k_{ij} \choose 2}  (n_i-1)^2 \]
 \begin{eqnarray*} \delta_j^3 &=&    (z_j  - 1)   
    - \sum_{i=1}^r  \left( \frac{p_{ij}}{d_{ij}} {k_{ij} \choose 2}   (n_i^{a_i} - 1)(n_i-1)\right. \\
    & &{}+\left( {k_{ij} \choose 3}  -  \frac{p_{ij}}{d_{ij}} {k_{ij} \choose 2} {a_i \choose 2}\right) (n_i-1)^3  \bigg)  + \delta_j^4 
     \end{eqnarray*}
for some $\delta_j^4 \in \Lambda_4$. With these values of $\delta_j^2,\delta_j^3$ the formula for $\xi_3$ in Theorem \ref{KD4/KD5} turns into the one we wished to prove.\hfbox\V

\section{Proofs for section 3}

All results quoted in section 3 are based on the following result.

\begin{satz}\label{Torkrit}\quad Let $\Delta$ be a descending filtration of $I(N)$ by subgroups 
$\IZ(N) = \Delta_1
\supset \Delta_2 \supset \ldots$ such that ${\rm Tor}_1^{\mathbb{Z}}\Big( 
\frac{\dst \Delta_1}{\dst \raisebox{-0.6mm}{$\Delta_{n-i}$}}\,,\,
\frac{\dst \IZ^i(T)}{\dst \IZ^{i+1}(T)} \Big) = 0$ for $1 \le i \le n-2$. 
Then there is a
natural isomorphism
  \[ \frac{\dst {\cal K}_n^{\Delta}}{\dst {\cal K}_{n+1}^{\Delta}} \hspace{2mm}\cong\hspace{2mm}
\bigoplus_{i=1}^{n-1} \hspace{2mm}
 \left( \frac{\dst \Delta_{n-i}}{\dst \Delta_{n-i+1}} \right) \hspace{1mm}\ot\hspace{1mm} \left( \frac{\dst \IZ^i(T)}{\dst
\IZ^{i+1}(T)}\right)  \]
induced by multiplication in $\IZ(G)$ (from the right to the left).
\end{satz}

\proof Let $1 \le i \le n-1$ and consider the following diagram where we use the notation of Proposition \ref{amalgam}, and where $\hat{\mu}_i$ is the factorization of $(\iota_1 {\small  \circ }  \cdots {\small  \circ } 
\iota_{i-1})\mu_i$
through its image.
  \[ \begin{matrix} 
\Imm{\Delta_{n-i}\IZ^i(T)}  & \stackrel{}{\hra} & 
\frac{\dst \Delta_1 \IZ (T)}{\dst \sum_{k=1}^n \Delta_{n-k+1}\,\IZ^{k}(T)} & & \cr
\surup{\hat{\mu}_i}  &  &  \mapup{\iota_1 {\small  \circ }  \cdots {\small  \circ } 
\iota_{i-1}} & &    \cr
 \frac{\dst \ruleu\Delta_{n-i}}{\dst \Delta_{n-i+1}}  \hspace{1mm}\ot\hspace{1mm} \frac{\dst \IZ^i(T)}{\dst
\IZ^{i+1}(T)}   &  \mr{\mu_i}  & 
\frac{\dst \Delta_1 \IZ^{i }(T)}{\dst \sum_{k=i}^n \Delta_{n-k+1}\,\IZ^{k}(T)}   & 
\ml{\iota_i}  & \frac{\dst \Delta_1 \IZ^{i+1}(T)}{\dst \sum_{k=i+1}^n
\Delta_{n-k+1}\,\IZ^{k}(T)} 
\end{matrix}\]
If $i=1$, skip the upper part of the diagram. Note that the Tor-term in the hypothesis also
vanishes for $i=n-1$ since then $\Delta_1/\Delta_1 = 0$. Then Proposition  
\ref{amalgam} and the hypothesis imply that  
   \BE\label{stufei} \mbox{$\mu_i$ and $\iota_i$ are injective and  
$\Imm{\mu_i} \cap \Imm{\iota_i} = 0$.} \EE
 So $ (\iota_1 {\small \circ} \cdots {\small \circ} \iota_{i-1})$ is injective, and
consequently the map  $\hat{\mu}_i$ is an isomorphism. It remains to
show that in $\Delta_1\IZ(T) / {\cal K}_{n+1}^{\Delta}$,
  \[ \Imm{\Delta_{n-j-1}\IZ^{j+1}(T)}  \cap  \sum_{l=1}^j  \Imm{\Delta_{n-l}\IZ^l(T)} = 0 \]
for all $1 \le j \le n-2$. We actually prove more, namely that for all $1 \le j \le n-2$,
  \BE\label{ind} \Imm{\Delta_1 \IZ^{j+1}(T)}  \cap  \sum_{l=1}^j  \Imm{\Delta_{n-l}\IZ^l(T)} = 0
\,.\EE
We proceed by induction on $j$. For $j=1$, $\Imm{\Delta_1 \IZ^{2}(T)}  \cap    \Imm{\Delta_{n-1}\IZ(T)} = 
\Imm{\iota_1} \cap \Imm{\mu_1} = 0$ by \REF{stufei}. Now suppose that relation \REF{ind} is true for $j=i-1$.
 Let $x \in \Delta_1 \IZ^{i+1}(T)$
such that the coset $x + {\cal K}_{n+1}^{\Delta}$ lies in $\sum_{l=1}^i  \Imm{\Delta_{n-l}\IZ^l(T)}$. Then there are $y
\in 
 \frac{\dst \Delta_{n-i}}{\dst \Delta_{n-i+1}}  \hspace{1mm}\ot\hspace{1mm} \frac{\dst \IZ^i(T)}{\dst
\IZ^{i+1}(T)}$, $z \in \sum_{l=1}^{i-1}  \Imm{\Delta_{n-l}\IZ^l(T)}$ such that 
 \[  (\iota_1 {\small \circ} \cdots {\small \circ} \iota_{i })
(\bar{x}) =  x + {\cal K}_{n+1}^{\Delta}  =   (\iota_1 {\small \circ} \cdots {\small \circ} \iota_{i-1}) (\mu_i y) + z\,.\]
Thus 
\begin{eqnarray*}
(\iota_1 {\small \circ} \cdots {\small \circ} \iota_{i-1})( \iota_{i } \bar{x} - \mu_i) (y) ) &
\in &  \Imm{\Delta_1 \IZ^{i}(T)}  \cap  \sum_{l=1}^{i-1}  \Imm{\Delta_{n-l}\IZ^l(T)} \\
 &=& 0
\end{eqnarray*}
by the induction hypothesis. Therefore, $\iota_{i } (\bar{x}) = \mu_i(y)$ by injectivity of 
$\iota_1 {\small \circ} \cdots {\small \circ} \iota_{i-1}$, and whence $\bar{x} = 0$ by
\REF{stufei}. Thus  $ x + {\cal K}_{n+1}^{\Delta}  =  (\iota_1 {\small \circ} \cdots 
{\small \circ} 
\iota_{i })(\bar{x}) = 0$, so \REF{ind} also holds for $j=i$, as was to be shown.  \hfbox\V

   \comment{
Next we  exhibit conditions on $N$ and $T$ which ensure that the hypothesis of
Corollary \ref{Torkrit} holds. For this purpose, we need to generalize \cite[Theorem 1.7]{PolProp} to
arbitrary N-series, as follows.
}

Now we are ready to prove the results stated in section 3. Starting out from \REF{AugTallg} and \REF{FoxTallg} Theorem
\ref{Aug+FoxT} is an immediate consequence of Theorem \ref{Torkrit} whose hypothesis is satisfied here by Corollary
\ref{UGrbism}: in fact, for $1\le i\le k$ the group $I^i(T)/I^{i+1}(T)$ is torsion-free, while for $k+1\le i\le n-2$ we have $n-i\le n-k-1$, whence $\Delta_1/\Delta_{n-i}$ is torsion-free.

Corollary \ref{Aug+FoxTU} then follows, again using Corollary
\ref{UGrbism}.

 In order to prove Theorem \ref{FoxN} first use \REF{FoxNallg} to reduce to the computation of $\Gamma_{n-1}^*I(N)/\Gamma_{n}^*I(N)$. Now the conjugation isomorphism $I(N)/\Lambda_{n-i-1}I(N)\hcong I(N)/I(N)\Lambda_{n-i-1}$ (see
  Remark \ref{conjonFox}) shows that Theorem \ref{Torkrit} applies for $\Delta=I(N)\Lambda$; it provides the desired decomposition of 
 $\Gamma_{n-1}^*I(N)/\Gamma_{n}^*I(N)$ by means of the mirror symmetry device in Remark \ref{symbem}.
 
 Finally,  Corollary \ref{FoxNU} then follows, once again using 
Corollary \ref{UGrbism}.\hfbox


%


%


\begin{thebibliography}{12}

\bibitem {Ba-GrI}  F.\ Bachmann und
L.\ Gruenenfelder,  \"Uber Lie-Ringe von Gruppen und ihre universellen
Enveloppen,   \textit{Comment. Math. Helv.}  {\bf 47} (1972),  332-340.


\bibitem {Ba-GrII}  F.\ Bachmann und
L.\ Gruenenfelder,  Homological methods and the third dimension
subgroup,   \textit{Comment. Math. Helv.}  {\bf 47 } (1972), 526-531.

\bibitem {PBW} P.-P. Grivel, Une histoire du th\'eor\`eme de Poincar\'e-Birkhoff-Witt, \textit{Expo.\ Math.} {\bf 22} (2004), 145-184.


\bibitem {Q3}  M.\ Hartl,  Some successive quotients of group ring filtrations
induced by N-series,   \textit{Comm. in Algebra}  {\bf 23}  (10) (1995), 3831-3853.

\bibitem{PolProp}   M. Hartl,  Polynomiality properties of group extensions with
torsion-free abelian kernel,    \textit{J.\ of Algebra}  {\bf 179} (1996), 380-415.

\bibitem{GoG}     M. Hartl,     The nonabelian tensor square and Schur multiplicator of
nilpotent groups of class 2,   \textit{J.\ of Algebra}  {\bf 179}  (1996)  416-440. 

\bibitem {D3F2}  M.\ Hartl,  The relative second Fox and third dimension subgroup of arbitrary groups,    \textit{Indian J. Pure Appl. Math.} {\bf 39} (5) (2008), 435-451.

\bibitem {HMP}  M.\ Hartl, R.\ Mikhailov, I.\ B.\ S.\ Passi, Dimension quotients, \textit{J. Indian Math. Soc.} (N.S.) 2007, Special volume on the occasion of the centenary year of IMS (1907-2007), 63-107 (2008).

\bibitem{Fox2}  M.\ Hartl,  On Fox quotients of arbitrary group algebras, \textit{Internat.\ J. Algebra Comput.,} to appear; http://arxiv.org/abs/0707.0281.

 


\bibitem {Hartley}  B.\ Hartley , The residual nilpotence of wreath products,   
\textit{Proc.\ London Math.\ Soc.}  (3) {\bf 20} (1970), 365-392.

\bibitem {KV27}  R.\ Karan, L.R.\ Vermani,    A note on augmentation quotients,  
\textit{Bull.\ London Math.\ Soc.}   {\bf 18}  (1986)  5-6. 

\bibitem {KV29}  R.\ Karan, L.R.\ Vermani,   Augmentation quotients of integral
group rings,   \textit{J.\ Indian Math.\ Soc.}    {\bf 54 }(1989)   107-120. 

\bibitem {KV30}  R.\ Karan, L.R.\ Vermani,    Augmentation quotients of integral
group rings--II,   \textit{J.\ Pure Appl.\ Algebra}   {\bf  65}  (1990)   253-262. 

\bibitem {KV31}  R.\ Karan, L.R.\ Vermani,    Augmentation quotients of integral
group rings--III,    \textit{J.\ Indian. Math. Soc.} {\bf 58} (1992), 19-32.

\bibitem {KV32}  R.\ Karan, L.R.\ Vermani,  Corrigendum:  Augmentation quotients of integral
group rings--II,    \textit{J.\ Pure Appl.\ Algebra}    {\bf 77}  (1992),   229-230. 


\bibitem {Kh12}  M.\ Khambadkone,  Subgroup ideals in group rings I,    \textit{J.\ Pure
Appl.\ Algebra}    {\bf 30}  (1983)  261-275. 

\bibitem {Kh13}  M.\ Khambadkone,  On the structure of augmentation ideals in group
rings,    \textit{J.\ Pure Appl.\ Algebra}   {\bf 35}  (1985)  35-45. 

\bibitem {Kh14}  M.\ Khambadkone,  Augmentation quotients of semidirect products,  
\textit{Arch.\ Math.}   {\bf 45}  (1985)  407-420. 

\bibitem {Kh15}  M.\ Khambadkone, Subgroup ideals in group rings III, \textit{Communications in Algebra} {\bf 14} (3) (1986), 411-421.

\bibitem{Lo} G.\ Losey, On the structure of $Q_2(G)$ for finitely generated groups, \textit{Canad.\ J.\ Math.} {\bf 25} (1973), 353-359.

\bibitem{ML}  S.\ Mac Lane,   \textit{Homology,} Springer Grundlehren, Vol.
114, Springer-Verlag Berlin-G\"ottingen-Heidelberg, 1963.

\bibitem{MP}  R.\ Mikhailov and  I.\ B.\ S.\ Passi, 
Augmentation powers and group homology,  \textit{J.\ Pure Appl.\ Algebra}  {\bf 192} (2004), 225-238.
 
 \bibitem{MPbook}  R.\ Mikhailov and  I.\ B.\ S.\ Passi, 
\textit{Lower central and dimension series of groups,} Lecture Notes in Mathematics, Vol.\ 1952, Springer Verlag Berlin-Heidelberg, 2009.

\bibitem{PaPF}    I.B.S. Passi,  Polynomial functors, \textit{Proc.\ Cambridge Philos.\ Soc.} {\bf 66} (1969), 505-512.


\bibitem{Pa}    I.B.S. Passi,   \textit{Group Rings and Their Augmentation
Ideals,}   Lecture Notes in Math., Vol. 715, Springer-Verlag, Berlin,
Heidelberg, New York, 1979.


\bibitem{Qu}  Quillen, D.,    On the associated graded of
a group ring, \textit{J. Algebra} {\bf 10} (1968),  411-418.

\bibitem{Rotman}  J.\ J.\ Rotman,  {\em   An Introduction to Homological Algebra}\/,  Academic
Press, Inc., 1979.


\bibitem{SaR}   R.\ Sandling, Dimension subgroups over arbitrary coefficient rings, \textit{J. Algebra} {\bf 21} (1972), 250-265. 

\bibitem{Sa}   R.\ Sandling, The dimension subgroup problem, \textit{J. Algebra} {\bf 2} (1972), 216-231. 

\bibitem{Sa-Ta}   R.\ Sandling and K.-I.\ Tahara,  Augmentation quotients of
group rings and symmetric powers, {\em Math.\ Proc.\ Camb.\ Phil.\ Soc.} {\bf 85} (1979), 
247-252. 

\bibitem{Ta3}    K.-I.\ Tahara,  On the structure of $Q_3(G)$ and the fourth dimension subgroups, {\em Japan J. Math.} {\bf 3} (1977),  381-394.

\bibitem{Ta4}    K.-I.\ Tahara, The augmentation quotients of group rings and the fifth dimension subgroups, \textit{J.\ Algebra} {\bf 71} (1) (1981), 141-173.

\bibitem{Ta}    K.-I.\ Tahara,  Augmentation quotients and dimension subgroups of
semidirect products, {\em Math.\ Proc.\ Camb.\ Phil.\ Soc.} {\bf 91} (1982),  39-49.


\bibitem {KV47} L.R.\ Vermani and  R.\ Karan,   Augmentation quotients of integral
group rings, III-- Corrigendum,  \textit{J.\ Indian. Math. Soc.} {\bf 59} (1993), 261-262.

\bibitem {Ve} L.R.\ Vermani, Augmentation quotients of integral group rings, in \textit{Groups--Korea '94 (Pusan),} 303-315, de Gruyter, Berlin, 1995.

\bibitem {Whitcomb}  A.\ Whitcomb,  The group ring problem, \textit{Thesis}, University of
Chicago, March 1968.
\end{thebibliography}
\end{document}